\documentclass{amsart}

\usepackage{amssymb}

\title{Finite symmetry groups in physics}
\author{Robert Arnott Wilson}
\date{15th February 2021; revised 28th July 2021; this version 8th November 2023}
\address{Queen Mary University of London}
\email{r.a.wilson@qmul.ac.uk}

\newcommand{\RR}{\mathbf R}
\newcommand{\CC}{\mathbf C}
\newcommand{\HH}{\mathbf H}
\newcommand{\ZZ}{\mathbf Z}
\newcommand{\QQ}{\mathbf Q}

\begin{document}
\begin{abstract}
Finite symmetries abound in particle physics, from the weak doublets and generation triplets to
the baryon octet and many others. These are usually studied by starting from a Lie group, and
breaking the symmetry by choosing a particular copy of the Weyl group. I investigate the
possibility of instead taking the finite symmetries as fundamental, and building the Lie groups
from them by means of a group algebra construction.

Finite group algebras are the natural algebraic structures in which finite symmetry groups, such as 
the symmetry group of three generations of elementary fermions, are embedded in Lie groups, that are
necessary for the formalism of quantum field theory, including the gauge groups of the fundamental forces. 
They are also the natural algebraic structures for describing representations of groups, which are
used for describing elementary particles and their quantum properties.
It is natural therefore to ask the question whether finite group algebras can provide a formal underpinning
for the standard model of particle physics, and if so, whether this foundation can explain any aspects of the
model that are otherwise unexplained, such as the curious structure of the combined gauge group,
or the mixing angles between the different forces.

In this paper I investigate the relationships between finite symmetry groups and the gauge groups
of each of the fundamental forces individually and in combination,  
and show
that
 the geometry of representations of the finite groups can be used to
predict accurate values for a number of the mixing angles in the standard model, including the
electro-weak mixing angle, one lepton mixing angle, one quark mixing angle and one 
of the CP-violating phases.

\end{abstract}
\maketitle

\section{Introduction}
\label{intro}
\subsection{Historical motivation}
\label{history}
Einstein made a number of very pertinent remarks about the foundations of quantum mechanics, often summarised in the famous
 phrase ``God
does not play dice". For example, in 1935 he wrote \cite{Einstein1935}:
\begin{quote}
``
In any case one does not have the right today to maintain that the foundation must consist in a field theory in the sense of Maxwell.
The other possibility, however, leads in my opinion to a renunciation of the time-space continuum and a purely algebraic physics.
Logically this is quite possible [...] Such a theory doesn't have to be based on the probability concept."
\end{quote}
In 1954, he went further, and opined \cite{Einstein1954}:
\begin{quote}
``I consider it quite possible that physics cannot be based on the field concept, i.e., continuous structure.
In that case, nothing remains of my entire castle in the air, gravitation theory included."
\end{quote}
Nevertheless, modern physics has stuck with the field concept, and his `castle in the air' remains in general use.

In 1935, there simply was not enough experimental evidence to provide much of a clue as to the possible structure
of a `purely algebraic physics', beyond the obvious fact that it must contain a copy of the quaternion group. But today, there is
an enormous amount of experimental evidence, that essentially shows there is (almost) a unique possibility for a purely
algebraic physics \cite{perspective}. Whether it actually works or not, is a separate question. But we should at least try it.

\subsection{Paradigms}
One of the deepest problems in the theory of elementary particles is to explain why there are three generations
of elementary fermions. The standard model is fundamentally a massless theory of a single generation, and the
Higgs mechanism for providing mass does not really explain why there are three different masses for the electron in particular,
and certainly does not explain the values of those masses.
Since it does not seem to be possible for the standard model to explain the existence of three generations, why not turn the question around,
and ask if the existence of three generations can explain the standard model?

The point of asking the question in this form is that there is no real explanation for why the gauge group of the standard model is what it is,
other than the fact that it works, and makes accurate predictions. Many attempts have been made over the past half-century
to derive this gauge group from something else, starting with the Grand Unified Theory (GUT) of Georgi and Glashow \cite{GG}
and that of  Pati and Salam \cite{PatiSalam} in 1974. In general, such attempts usually start from a larger gauge group or algebra,
and seek to explain the standard model via a process of `spontaneous symmetry-breaking'. 

The GUT paradigm supposes that this
larger structure is a Lie group or Lie algebra \cite{MDW,Chester}, but other models using Clifford algebras \cite{Lu,Stoica,Clifford} 
have also been popular, not to mention
more exotic algebras including Jordan algebras \cite{MD,Todorov}, division algebras \cite{GuGu,GDixon,CFurey,NFurey}
and so on. Despite some successes, such
models have generally either failed to be predictive, or have predicted things (such as proton decay, or new particles) 
that have not been detected. An alternative approach is to try to build up to the standard model from something smaller,
such as a finite symmetry group,
rather than restrict down from something bigger.

\subsection{Group algebras and Clifford algebras}
This approach can be considered a generalisation of the Clifford algebra approach, since Clifford algebras are the special case in which
the finite symmetry group is a double cover of a direct product of some number of copes of the group $Z_2$ of order $2$.
In the language of finite group theory, in even dimensions, these groups are called \emph{extraspecial $2$-groups}, and in two dimensions
are just the quaternion group $Q_8$ and the dihedral group $D_8$, both of order $8$. The Clifford algebra arises as the faithful 
(i.e. `spinor') part
of what turns out to be the \emph{real group algebra} of the finite group. 

In the case of $Q_8$, one must choose two group elements of order $4$ as generators, squaring to $-1$,
so that the corresponding Clifford algebra is $Cl(0,2)$, isomorphic to the quaternion algebra $\HH$. 
In the case of $D_8$, one can choose either two elements of order $2$, or one of order
$2$ and one of order $4$, so that both $Cl(2,0)$ and $Cl(1,1)$ arise here, and both are isomorphic to the $2\times 2$ real matrix algebra
$M_2(\RR)$. In terms of Pauli matrices, the three possible choices of generators are 
\begin{itemize}
\item
$(\sigma_1,\sigma_2)$ 
for $Cl(2,0)$, 
\item $(\sigma_1,i\sigma_2)$
for $Cl(1,1)$, and 
\item $(i\sigma_1,i\sigma_2)$ 
for $Cl(0,2)$.
\end{itemize}

In four dimensions, the two algebras are $M_4(\RR)$ and $M_2(\HH)$. The former arises from the group 
\begin{align}
(D_8 \times D_8)/Z_2 = (Q_8 \times Q_8)/Z_2
\end{align}
in which different choices of generators give rise to $Cl(3,1)$ and $Cl(2,2)$.
Similarly, the latter arises from the group
\begin{align}
D_8\otimes Q_8 : = (D_8\times Q_8)/Z_2
\end{align}
in which various choices of generators give rise to $Cl(1,3)$, $Cl(4,0)$ and $Cl(0,4)$.
In terms of Dirac matrices, one possible choice of generators is as follows:
\begin{itemize}
\item $\gamma_0,i\gamma_1,i\gamma_2,i\gamma_3$ for $Cl(4,0)$,
\item $i\gamma_0,i\gamma_1,i\gamma_2,i\gamma_3$ for $Cl(3,1)$,
\item $\gamma_0,\gamma_1,i\gamma_2,\gamma_3$ for $Cl(2,2)$,
\item $\gamma_0,\gamma_1,\gamma_2,\gamma_3$ for $Cl(1,3)$,
\item $i\gamma_0,\gamma_1,\gamma_2,\gamma_3$ for $Cl(0,4)$
\end{itemize}

In $2n+1$ dimensions, one adds another generator that commutes with the first $2n$ generators, and either squares to $+1$ or to $-1$.
The former case gives a direct product of groups, and two copies of the smaller algebra, while the latter gives a so-called
`central' (or tensor) product, and a complexification of the smaller algebra. For example, in three dimensions the cases are
\begin{itemize}
\item $Z_2 \times D_8$, giving rise to $M_2(\RR)+M_2(\RR)$ and $Cl(2,1)$,
\item $Z_2 \times Q_8$, giving rise to $\HH+\HH$ and $Cl(0,3)$,
\item $Z_4\otimes D_8 = Z_4\otimes Q_8$, giving rise to $M_2(\CC)$ and both $Cl(3,0)$ and $Cl(1,2)$. 
\end{itemize}
In five dimensions one obtains the Dirac algebra (complex Clifford algebra) $M_4(\CC)$ from the group
\begin{align}
Z_4\otimes D_8\otimes D_8 = Z_4\otimes D_8\otimes Q_8 = Z_4 \otimes Q_8\otimes Q_8,
\end{align}
in which there is a wide choice of generators for $Cl(4,1)$, $Cl(2,3)$ and $Cl(0,5)$, from which one can construct
either de Sitter, anti-de Sitter or Euclidean spacetimes according to preference. The reversed signatures
$Cl(5,0)$ and $Cl(1,4)$ arise from $Z_2 \times D_8\otimes Q_8$, while $Cl(3,2)$ arises from $Z_2 \times Q_8\otimes Q_8$.

For the description of triplets, such as generations or colours, Clifford algebras are insufficient, since there are no natural
triplet symmetries in Clifford algebras. The idea of the present paper, therefore, is to add a triplet symmetry to the finite group
that defines the Clifford algebra of choice, and to see where that leads. The obvious choice is to use the group
$Z_4\otimes Q_8\otimes Q_8$ of order $64$, 
which has two commuting triplet symmetries, which could implement colour and generation symmetries
respectively, in a group algebra in $576$ dimensions. 

For a first toy model, however, I take instead the smallest case which has a non-trivial
triplet symmetry, that is $Q_8$, giving rise to a group of order $24$ and a $24$-dimensional algebra.
It is then possible to obtain an algebra of dimension $24\times 24=576$ by treating left and right multiplications in the algebra
independently. 
Such an algebra therefore has enough degrees of freedom to add both 
colours and generations to the Dirac algebra, and hence to unify the entire Standard Model in a single algebraic construct,
that is a natural generalisation of a Clifford algebra.

\subsection{Plan of the paper}
Since finite group algebras are not familiar constructions in particle physics, I begin by studying some toy models in order to
illustrate the methods, and to demonstrate that the approach has the potential to explain certain things that other approaches cannot explain.
The smallest possible example is the group of order $3$, studied in Section~\ref{Z3model}, with potential applications to
both the triplet colour symmetry of quantum chromodynamics (QCD) and the triplet generation symmetry of fundamental fermions.
This group is a discrete version of a $U(1)$ gauge group, which is in turn a subgroup of the finite group algebra.
The second example, in Section~\ref{Z4symmetry}, 
is the cyclic group of order $4$, which is the only other discrete version of $U(1)$ whose group algebra contains just one copy of
$U(1)$.

This is followed by an extended discussion of the quaternion group $Q_8$ as a finite version of $SU(2)$. 
This group can be used (Section~\ref{Q8spin})
to study both
the (non-relativistic) spin of elementary particles, and 
the weak interaction, both of which are described in the Standard Model by the group $SU(2)$. 
This is followed (Section~\ref{massq}) by a description of
how the group algebra can be used to quantise certain masses, and derive three mixing angles.

We then put $Z_3$ and $Q_8$ together, to form the binary tetrahedral group of order $24$, which is the subject of
Sections~\ref{tetrahedron} and \ref{EWunif}. 
Its group algebra is therefore a $24$-dimensional real algebra, made up of analogues of 
\begin{itemize}
\item 
the $12$-dimensional
gauge group $(U(1)\times SU(2)\times SU(3))/Z_6$ of the Standard Model, 
\item 
the $6$-dimensional Lorentz group
$SL(2,\CC)$ and 
\item 
four real scalars and one complex scalar. 
\end{itemize}
The algebra therefore contains $24$ apparently 
arbitrary parameters, which I shall attempt to elucidate to some extent.
In Section~\ref{SM} I show how to derive the main structures of the Standard Model from this group algebra, 
especially the strong force and the Dirac equation.

Section~\ref{physics} is a general discussion of the emergence of classical physics from the discrete interactions that are implicit in the
model.  
This includes the emergence of spacetime itself, and consequences for quantum theories of gravity.
In Section~\ref{waveparticle} I consider various aspects of quantum field theory, 
including wave-particle duality and the allocation of particles to different parts of the model.
Finally in Section~\ref{other} I discuss relationships to other models.

\section{A triplet symmetry toy model}
\label{Z3model}
\subsection{The group and its algebra}
Let us therefore begin with the assumption that there are three generations of electrons, and use the finite group $Z_3$ of order $3$
to study the symmetry \cite{Zee}. There are three elements $e,f,g$ in the group, satisfying the multiplication rules
\begin{align}
ee=fg=gf=e,\cr ff=eg=ge=g,\cr gg=ef=fe=f.
\end{align}
To get a continuous theory (for example, a quantum field theory \cite{WoitQFT}) out of this finite group
we must  
take (at least) all linear combinations
\begin{align}
xe+yf+zg
\end{align}
with $x,y,z$ real numbers. 

Two such linear combinations can be multiplied together by combining the group multiplication with ordinary
arithmetic to obtain
\begin{align}
(xe+yf+zg)(ue+vf+wg) = & (xu+yw+zv)e + \cr
 & (xv+yu+zw)f + (xw+yv+zu)g.
\end{align}
Hence we obtain a $3$-dimensional algebra, called the real group algebra 
of $Z_3$, and denoted $\RR Z_3$, which happens to be isomorphic (as an algebra) to the sum $\RR+\CC$ of the real and complex numbers
\cite[Chapter 6]{JamesLiebeck}.
To see this isomorphism explicitly, we can take $e+f+g$ as a generator for the real subalgebra $\RR$, noting that
\begin{align}
(e+f+g)^2 = 3e+3f+3g,
\end{align}
so that $(e+f+g)/3$ acts as the real number $1$. 

Then we can take $2e-f-g$ and $f-g$ as generators for $\CC$, and compute
\begin{align}
(f-g)^2 & = f+g-2e\cr
(2e-f-g)^2 & = 6e-3f-3g\cr
(f-g)(2e-f-g)  & = 3f-3g 
\end{align}
so that $(2e-f-g)/3$ acts as the complex number $1$, and $f-g$ is a square root of $-3$. Finally we compute
\begin{align}
(e+f+g)(f-g) & = 0\cr
(e+f+g)(2e-f-g) & = 0
\end{align} 
to confirm that the algebra really is a direct sum of $\RR$ and $\CC$.

\subsection{Physical interpretations}
\label{toyelectrons}
If we now take out the two real scale factors from the algebra, what is left is simply the group $U(1)$ of unit complex
numbers, which we can interpret as one of the gauge groups of the standard model,
either $U(1)_{em}$ of (quantum) electrodynamics or $U(1)_Y$ of weak hypercharge, or perhaps something else, according to context. 

The subalgebra $\RR$ generated by $e+f+g$
is the fixed space of the triplet symmetry, so represents the common feature of the three generations---that is, let us say, the charge. 
This charge is not affected by the action of $U(1)$. 
The subalgebra $\CC$, on the other hand, represents the differences between the three generations---that is, essentially, the mass.
Hence the group $U(1)$ can in principle affect the mass, depending on how exactly we map from particles into $\CC$.

If, for example, we project onto the subalgebra $\CC$,  
removing the charge from the model, then we are left with a neutral particle, most
reasonably interpreted in this context as a neutron.  
We obtain $(2e-f-g)/3$ as the real part of the complex number,
and if we split this up into the three terms $2e/3$, $-f/3$ and $-g/3$ then we can compute the charges of these components as
\begin{align}
2e/3.(e+f+g) &= 2/3(e+f+g)\cr
-f/3.(e+f+g) & = -1/3(e+f+g)\cr
-g/3(e+f+g) & = -1/3(e+f+g),
\end{align}
where $e+f+g$ is (as already noted above) the  
measurement scale  
for charge.
Thus these three components correspond to, or can be interpreted as, 
the up quark and two down quarks that make up the neutron. 

One might naively expect that $e,f,g$ themselves (or more precisely, their negatives)
represent the three generations of electron, but this appears not to be the case. The most
reasonable interpretation seems to be that $e-f-g$ represents the ordinary electron, and $-e+f-g$ and $-e-f+g$ represent the other two generations,
that is the muon and the tau. 
Since this toy model makes no mention of neutrinos, or of baryonic mass, but only charge, it makes sense to subtract an electron $e-f-g$ from a neutron
$(2e-f-g)/3$ to get a proton $(-e+2f+2g)/3$, consisting of one down quark (charge $-1/3$) and two up quarks (charge $2/3$) as in the standard model.
Thereby $e$, $f$ and $g$ acquire an interpretation as quarks,
or, 
more abstractly, as the
\emph{colours} of the three quarks. 
Colour confinement can then be said to arise from the splitting of the colour algebra $\CC$ from the charge algebra $\RR$,
or directly from the group-theoretic property that $e=fg$. The electrons are then colourless, since they contain the three colours $e,f,g$ 
in equal measure (ignoring the signs, which, as we have seen, negate electric charge but do not change the colour).

\subsection{A mixing angle}
\label{massplane}
Now suppose that we use the algebra instead to represent three generations of electrons.
We first project out the charge on the electrons, so that they are represented by the vertices of an equilateral triangle in a real $2$-space.
This real $2$-space
no longer has a natural identification 
with the complex numbers, since the particles lie in representation spaces of the algebra, not in the algebra itself. The differences
between the three generations lie in the mass alone, so that this $2$-space is a `mass plane', in which the usual scalar (inertial) mass can be
measured in a specific direction. 
Figure~\ref{emufig} shows the mass axis drawn through the $\tau$ vertex, with projections from $e$ and $\mu$ onto this line, in order
to compare the mass differences $m(\tau)-m(\mu)$ and $m(\tau)-m(e)$. 
\begin{figure}
\caption{The electron-muon mixing angle\label{emufig}}
\begin{picture}(200,220)
\put(0,20)
{
\put(20,90){\circle*{5}}
\put(170,0){\circle*{5}}
\put(170,180){\circle*{5}}
\put(170,0){\line(0,1){180}}
\put(20,90){\line(5,3){150}}
\put(20,90){\line(5,-3){150}}
\put(70,150){\vector(2,-3){115}}
\put(20,90){\line(3,2){70}}
\put(170,180){\line(-3,-2){100}}
\put(10,90){$e$}
\put(175,180){$\mu$}
\put(175,0){$\tau$}
\put(160,20){$\theta$}
\put(105,40){$60^\circ-\theta$}
\put(182,-15){mass}
\put(80,137){$0$}
}
\end{picture}
\end{figure}
The diagram is drawn approximately to scale, 
within the limitations of the \LaTeX\ picture environment,
which mean that the angle at the $\tau$ vertex, marked as $60^\circ$, is just over $59^\circ$ in the picture.

The scale 
is approximately 30MeV$/c^2$ to a millimetre, so that
the mass of the electron is indistinguishable from $0$, as marked on the mass axis.
A simple calculation yields the formula for the angle $\theta$: 
\begin{align}
\frac{\cos(60^\circ-\theta)}{\cos\theta}  & = \frac{m(\tau)-m(e)}{m(\tau)-m(\mu)}\cr
\Rightarrow \frac{\sqrt{3}}{2}{\tan\theta} +\frac12 &\approx 1.06292\cr
\Rightarrow \tan\theta & \approx .65000\cr
\Rightarrow \theta & \approx 33.024^\circ.
\end{align}
We expect this angle, or something like it, to appear somewhere in the standard model, most likely as a
lepton-mixing angle between the first and second generations. In other words, it should be closely related to the
mixing angle $\theta_{12}$ in the Pontecorvo--Maki--Nakagawa--Sakata (PMNS) matrix \cite{Pontecorvo,MNS}.

The experimental value of this angle is
\begin{align}
\theta_{12} & = {33.41^\circ}^{+.75^\circ}_{-.72^\circ}
\end{align}
in complete agreement with the value obtained geometrically from the lepton masses.
Either this is an extraordinary coincidence, or it is an indication that the group algebra method is a powerful method for
making predictions in particle physics. (Incidentally, the closeness of $\tan\theta$ to $13/20$ is \emph{not} a prediction of this model,
and appears to be a pure coincidence. The approximation used in the diagram is $\tan\theta=2/3$, 
corresponding to an angle $\theta$ of approximately $33.69007^\circ$, which is, perhaps ironically, or perhaps not, even closer to
the experimental value than the value calculated above.)

\subsection{Commentary}
In the standard model, the parameter $\theta_{12}$ is a mixing angle between the electron neutrino and the muon neutrino,
and is supposed to be a property of neutrinos alone. In particular, the value of this parameter is, at least in principle,  
independent of
the masses of the three generations of electron. The group algebra approach, however, suggests that
these parameters may not in fact be independent.  
If so, then we obtain a more precise prediction of
$33.024^\circ$ for this particular neutrino mixing angle, which can possibly be tested by new neutrino-oscillation measurements.

The group algebra approach does not explain why it is the electron-muon mixing angle that appears here, and not the
electron-tau or muon-tau angles. But since the geometry is $2$-dimensional, only one of these parameters can appear.
In order to obtain the other two we require (at least) $3$-dimensional geometry. Another indication that we need a higher-dimensional space
is the fact that the zero of the mass scale does not coincide with the centre of symmetry. 

\section{From gauge groups to finite symmetries}
\label{Z4symmetry}
\subsection{Quadruplet symmetry}
We have seen that a triplet symmetry gives rise via the group algebra construction to a gauge group $U(1)$.
It is natural to ask if there are any other groups with this property? That is, which finite groups have a group algebra consisting of a single copy
of $U(1)$, together with some number of real scalars? It is actually quite easy to see that there is exactly one more, that is the cyclic group
$Z_4$ of order $4$. 

To prove this, note first that the group algebra is commutative, so the group is abelian, and is therefore a
direct product of cyclic groups. Any cyclic factor
of order greater than $4$ gives at least two copies of $U(1)$. Any cyclic factor of order $2$ doubles the number of copies of $U(1)$. Hence there is a single cyclic factor, which has order $3$ or $4$. The group algebra of $Z_4$ has the structure $\RR+\RR+\CC$, compared to $\RR+\CC$ for $Z_3$.
In addition to the trivial representation, there is a `doublet' representation, and the complex representation as rotation symmetries of a square.

One place where a cyclic symmetry of order $4$ might plausibly underlie the standard model is in the four types of fundamental fermion---that
is, the electron, neutrino, up and down quark. If we arrange these particles on a square in such a way that a reflection
symmetry of the square swaps each
weak doublet---that is the `left-handed' electron/neutrino and up/down quark doublets, with the negative charge on the left, 
say---then they appear in the cyclic order $\nu,e,d,u$, with charges $0,-1,-1/3,2/3$
respectively. 
This square is drawn in Fig.~\ref{Wfig}, with axes for the charge $Q$, weak hypercharge $Y_W$ and third component of weak isosopin $T_3$.

\begin{figure}
\caption{The electro-weak mixing angle\label{Wfig}}
\begin{picture}(200,220)
\put(0,120){\circle*{5}}
\put(80,0){\circle*{5}}
\put(120,200){\circle*{5}}
\put(200,80){\circle*{5}}
\put(0,120){\line(3,2){120}}
\put(0,120){\line(2,-3){80}}
\put(120,200){\line(2,-3){80}}
\put(50,175){\vector(2,-3){110}}
\put(80,0){\line(3,2){120}}
\put(215,120){$T_3$}
\put(45,20){\vector(3,2){170}}
\put(-20,70){\vector(1,0){250}}
\put(0,120){\line(0,-1){50}}
\put(80,0){\line(0,1){70}}
\put(200,80){\line(0,-1){10}}
\put(120,200){\line(0,-1){130}}
\put(-5,60){$-1$}\put(70,75){$-1/3$}
\put(190,60){$2/3$}
\put(-15,115){$e_L$}\put(65,-4){$d_L$}\put(205,75){$u_L$}
\put(107,198){$\nu$}
\put(215,60){$Q$}
\put(160,20){$Y_W$}
\put(135,50){$1/3$}\put(123,77){$0$}\put(30,30){$-1/2$}\put(175,120){$1/2$}
\put(40,160){$-1$}
\put(110,185){$\phi$}
\end{picture}
\end{figure}

\subsection{Calculating the mixing angle}
We can now calculate 
the angle $\phi$ between the charge axis and the weak hypercharge axis.
The calculation  
is rather easier than in the $Z_3$ case, since no angles of $60^\circ$ are involved, and reduces quickly to
\begin{align}
\tan\phi & = 3/2\cr
\Rightarrow \phi 
& \approx 56.30993^\circ
\end{align}
from which we find an approximation to the electro-weak mixing angle as
\begin{align}
\phi/2 & \approx 28.15497^\circ
\end{align}
compared to the CODATA 2018 \cite{CODATA2018} value
\begin{align}
\theta_W & = 28.172^\circ \pm .021^\circ.
\end{align}

The reason why we see the angle $2\theta_W$ rather than $\theta_W$ itself is that we are dealing here with fermions, not bosons.
Thus the projection from the charge axis onto the weak hypercharge axis gives rise to a factor of
\begin{align}
\cos \phi = 1 - 2 \sin^2(\phi/2)
\end{align}
which multiplies the charge in the vertex factors in the Feynman calculus for weak force neutral current interactions.
Conventionally, only the term $-2\sin^2(\phi/2)$ is used, and the constant is subsumed into other parts of the calculation.
This is related to the fact that the origin of coordinates is offset from the centre of symmetry.
The `neutral charge' itself must be taken perpendicular to the electric charge $Q$, and is represented in this diagram by the neutrino.

Reversing the above calculations gives 
\begin{align}
\sin^2\theta_W & = (1 -\cos\phi)/2\cr
& = 1/2 - 1/\sqrt{13}\cr
& \approx  .22265
\end{align}
compared to the CODATA 2018 
recommended value $.22290\pm.00030$.
We can also translate to the W/Z mass ratio
\begin{align}
\frac{m(W)}{m(Z)} & = \cos\theta_W\cr
& = \sqrt{1/2+1/\sqrt{13}}\cr
 & \approx .8816746
\end{align}
which is again in close agreement with both the standard model and experiment.
Note in particular that the group algebra approach does \emph{not} predict the so-called W/Z mass anomaly
that has been announced recently \cite{WZ}. 

\subsection{Taking stock}
The diagram is drawn accurately to scale, and makes clear that the 
origin of coordinates is not a natural zero for the neutral charge. 
A better convention might be that the neutral charge of the electron
is $0$, and is scaled so that the neutral charge of the down quark is $-1$. Then the up quark and neutrino have neutral charges
$-1/3$ and $2/3$ respectively, and the pattern of neutral charges is simply a rotation of the pattern of electric charges.
Thus the group $Z_4$ becomes a symmetry group uniting the electric and neutral charges.
The inverse symmetry applies if we take the neutrino to have neutral charge $-1$, the electron $-1/3$, the up quark $0$ and the down quark $2/3$.
It is not immediately obvious which is the more natural convention, or if either of them is consistent with experimental
measurements of neutral currents.

To summarise, we have found that the gauge group $U(1)$ arises from exactly two finite symmetry groups, $Z_3$ and $Z_4$, both of which
have natural interpretations in the standard model of particle physics, and both of which give elementary geometric derivations of
one of the otherwise unexplained mixing angles.
It may be worth noting that the angles $\theta$ in Fig.~\ref{emufig} and $\phi$ in Fig.~\ref{Wfig} add up to
$89.334^\circ$, which is experimentally (though not theoretically) consistent with being exactly $90^\circ$.
Probably this is a meaningless coincidence, but possibly not, if, for example, we find that both of these copies of $U(1)$
are identified with the same copy of $U(1)$ in the Standard Model. 
 It is time therefore to move on to the other gauge groups $SU(2)$ and $SU(3)$, and
see if they can also be related to finite groups.

\section{A finite version of $SU(2)$}
\label{Q8spin}
\subsection{The spin group and the quaternion group}
\label{quaternion}
Let us first ask whether there is a finite group whose group algebra consists of a single copy of $SU(2)$, together with some number of
copies of $\RR$. Such a group is necessarily non-abelian, and all its
cyclic quotients must be 
copies of $Z_2$.
Indeed, it is fairly easy to see that there is only one such group, namely the quaternion group $Q_8$ of order $8$. This group has elements
$\pm1$, $\pm i$, $\pm j$ and $\pm k$, satisfying the rules 
\begin{align}
&ij=k, jk=i, ki=j, \cr
&i^2=j^2=k^2=-1.
\end{align}

Indeed, serious group theory \cite{Zee}
first entered into quantum mechanics with the discovery of the spin of an electron, which has a direction relative
to the ambient space, and therefore requires a group $SU(2)$ for its description. On a macroscopic scale, the sum of the spins of a large number of 
elementary particles 
gives rise to
the phenomenon of magnetism. But it is experimentally impossible to measure the
direction of spin of an individual elementary particle. 
All one can do is choose a direction, and measure whether the spin is `up' or `down' in that direction. 
Indeed, even for silver atoms, the Stern--Gerlach experiment \cite{SternGerlach}
demonstrates that the direction of spin is
quantised.
Therefore we must assume that the elementary particles have,
intrinsically, only finite symmetries. The smallest finite group that is available to describe spin symmetries is the
quaternion group $Q_8$ described above.

Now if we take a large number of copies of $Q_8$, and add them together, we obtain \cite{grouprings} 
the integral group ring $\ZZ Q_8$, in which all the copies of $Q_8$ are identical, meaning that all electrons are identical, and so on. In macroscopic physics, the individual quanta disappear from view, and it
becomes reasonable to approximate the scaled copy of the integers $\mathbb \varepsilon \ZZ$ by the real numbers $\RR$, and study the
group algebra $\RR Q_8$ instead. There is, up to equivalence, only one faithful representation of $Q_8$, in the quaternion algebra $\HH$,
and there are four
$1$-dimensional real representations, of which one is trivial and three are representations of the three quotients $Q_8/Z_4\cong Z_2$.
Therefore the group algebra has a canonical (Wedderburn) decomposition as 
\begin{align}
\RR Q_8 &\cong  4\RR + \HH.
\end{align}
where $\HH$ denotes the quaternion algebra of all real linear combinations of $1,i,j,k$, under the same rules of multiplication as above.
Taking out the five real scale factors from this algebra, we are left with the group $SU(2)$,
which is exactly what we require in order to describe
magnetic spin in ordinary macroscopic space.

\subsection{Analysis of the group algebra}
\label{analysis}
In the representation theory of finite groups, the usual convention is for the finite group to act on the algebra by right multiplication, and the continuous groups to act by left-multiplication.  This allows us to treat both the discrete and continuous symmetries simultaneously. 
The macroscopic behaviour, on the other hand, requires the continuous group to act on both sides,
by conjugation, so that $SU(2)$ acts on $\HH$ as $SO(3)$, that is, fixing the real part, and rotating the imaginary part as a Euclidean
$3$-space. 
It is standard to interpret this $3$-space as real physical space, in the case when $SU(2)$ represents magnetic spin.
It would be reasonable, then, to interpret the real part of $\HH$ as a (non-relativistic) time. This allows the `spin' to be modelled 
as something
that takes place 
in space and time, rather than in isolation.

But the finite group has given us something in addition, namely four copies of the real numbers. These 
presumably represent four physical 
particles, at least three of which we can
detect magnetically. 
These
four particles 
are acted on by the finite group on the right, in three cases to change the sign, so that these three 
appear in a pair of spin states. They are invariant under conjugation by $SU(2)$, so under symmetries of
macroscopic space. But they are not invariant under the left action of the finite group. Since we started out modelling spin of electrons
and protons, these had better be two of the particles, and the most plausible %identification 
interpretation of the other two is as the neutron
and (electron) neutrino. The first three appear in two spin states, the last in only one, as we know also from the Wu experiment \cite{Wu}.
Thus this action describes the weak interaction, in the form of beta decay:
\begin{eqnarray}
e+p &\leftrightarrow& n+\nu.
\end{eqnarray}

Note that this is an action of the finite group, not of the spin group.
In the standard model \cite{Griffiths}, the weak interaction is described instead by a second copy of $SU(2)$, commuting with the spin group.
It is of course perfectly possible to define an action of $SU(2)$ on this $4$-dimensional real space, and have this action
commute with the spin $SU(2)$, and try to build a model of physics on top.
But this action is not compatible with 
the group algebra.
The standard model therefore puts the weak $SU(2)$ elsewhere, and doubles up the Weyl spinor into a Dirac spinor
so that there will be no conflict between the finite action and the continuous action. By doing so, however, it breaks the connection
between the weak interaction and the spin group, which then has to be put back in by hand, in the form of a mixing
between the weak force and electromagnetism. 

In the group algebra, the four $1$-dimensional representations 
contain the following elements:
\begin{eqnarray}
1a:&& 1+(-1)+i+(-i)+j+(-j)+k+(-k);\cr
1b: && 1+(-1)+i+(-i)-j-(-j)-k-(-k);\cr
1c:&& 1+(-1)-i-(-i)+j+(-j)-k-(-k);\cr
1d:&& 1+(-1)-i-(-i)-j-(-j)+k+(-k).
\end{eqnarray}
We then see that the elements of the group can be interpreted, not as particles, but as quantum numbers of the particles. Two quantum numbers are sufficient to distinguish these four particles, so that for example, if $\lambda_i$ and $\lambda_j$ denote the coefficients of $i$ and $j$ respectively, then we can take the charge to be $(\lambda_j-\lambda_i)/2$, and the weak isospin to be $\lambda_j/2$. We do not need to specify the coefficient of $k$, that is $\lambda_k=\lambda_i\lambda_j$.

In a similar way, the other half of the group algebra has a basis consisting of
\begin{eqnarray}
t&:=& 1-(-1),\cr
x&:=& i-(-i),\cr
y&:=& j-(-j),\cr
z&:=& k-(-k).
\end{eqnarray}
The names are chosen with the intention of suggesting a potential quantisation of spacetime, with $z$ representing the direction of spin, and the $x,y$ plane relating to the charge and isospin. 
The (speculative) idea, which I shall explore in more detail below, is that the combination of spin, charge and weak isospin for a (large) collection of interacting elementary particles 
may be sufficient to define the ambient physical spacetime.

This 
example illustrates the general principle of how a macroscopic force, in this case magnetism, emerges from individual quanta. In particular,
it illustrates how a macroscopic continuous symmetry arises from a microscopic discrete symmetry. It also illustrates how a breaking of the symmetry
at the microscopic level arises from an interaction between elementary particles, in this case the weak interaction,
and gives rise to a breaking of symmetry also at the
macroscopic level. The quaternion group $Q_8$ on its own is therefore sufficient to model a united magneto-weak force. It is not sufficient for
a complete electro-magneto-weak unification, however, for which we need a larger group.
This larger group must be capable of dealing with the fact that the electron comes in three different mass states, that is in three different
generations. 

\subsection{The automorphism group}
\label{auto}
The symmetries of this finite model of spin are described by the automorphism group of $Q_8$, that is isomorphic to
$Sym(4)$, the symmetric group on $4$ letters.
The automorphisms are of five types, where I denote automorphisms as permutations of the letters $W,X,Y,Z$: 
\begin{itemize}
\item the identity element, of order $1$.
\item $3$ even permutations of order $2$, such as $(W,X)(Y,Z)$; these are inner automorphisms, that is conjugation by $i$, $j$ and $k$. 
\item $8$ even permutations of order $3$, such as $(X,Y,Z)$; these can be represented as conjugation by unit quaternions $(-1\pm i\pm j\pm k)/2$.
\item $6$ odd permutations of order $2$, such as $(W,X)$; represented as conjugation by $i\pm j$, $j\pm k$ and $k\pm i$.
\item $6$ odd permutations of order $4$, such as $(W,X,Y,Z)$; represented as conjugation by $1\pm i$, $1\pm j$ and $1\pm k$.
\end{itemize} 

These symmetries may have many different interpretations, depending on what the original copy of $Q_8$ represents. The kind of symmetry that
it might be useful to look for is an extension of the neutrino/electron/proton/neutron symmetry discussed in the previous section, to
three generations of electron, plus the baryon octet \cite{GellMann}. 
In this way the baryons would exhibit a $3$-fold symmetry represented by the $3$-cycles,
that could perhaps be interpreted as a colour symmetry of the three constituent quarks \cite{QCD}. Similarly,
the electrons would divide into two transpositions, possibly representing the left-handed and right-handed spins. 

Alternatively, or in addition,
one might want to consider the odd permutations as representing leptons and quarks. The transpositions could then represent either
the left-handed and right-handed electrons, 
or the left-handed parts of both neutrinos and electrons. Similarly,
the $4$-cycles could represent quarks, either in up/down pairs 
or in left/right pairs.
Additional interpretations might arise from consideration of particle interactions, so that the identity element represents electromagnetic
interactions, 
 the $3$ elements of order $2$ represent the weak interaction,
and the $8$ elements of order $3$ the strong interaction. The strong interaction manifests itself both in massless form (as gluons) and
in massive form (as the meson octet, including pions and kaons).

Note also that the group $Z_4$ described in Section~\ref{Z4symmetry} acts as an automorphism group of this copy of $Q_8$,
and not as a subgroup. This is because $Z_4$ mixes quarks and leptons, while $Q_8$ does not. We might for example choose a generator for
$Z_4$ to be the automorphism that fixes $i$ and maps $j$ to $k$, so that it acts as a $4$-cycle on $(j,k,-j,-k)$. This generator squares to the
automorphism that negates $j$ and $k$, which is the same as the inner automorphism defined by conjugation by $i$. However,
it is not possible to identify this automorphism with the element $i$, since the former swaps quarks with leptons, and the latter does not.
Similarly, $Z_3$ also acts as an automorphism group of $Q_8$, for example as the $3$-cycle on $(i,j,k)$.
This leads in the Standard Model to two quite different actions of $U(1)$ on $SU(2)$, which are distinguished by using
$\gamma_5$ and the electro-weak mixing angle.

It may be worth mentioning in passing that if one extends the $S_4$ action by conjugation of the quaternions listed above,
to include both left and right multiplications, then one obtains the Weyl group of type $F_4$. This may be the fundamental reason
why the Lie group of type $F_4$, and various related types such as $D_4$ and $E_6$, are so attractive as potential ways
to extend the standard model \cite{E6,E8,DG,Furey,DMW}.

\subsection{Interpretation as a gauge group for the weak interaction}
As already noted,
one of the four copies 
of $\RR$ is the trivial representation of $Q_8$, while the other three are doublet representations (that is, representations of $Z_2$), 
which one can presumably interpret as weak doublets
for three separate generations of fermions. Just as in the case of the $Z_4$ doublet representation, these representations do not distinguish between
lepton doublets and quark doublets. Or perhaps it is better to say that they describe \emph{only} lepton doublets,
and quark doublets are somewhat different (as indeed the quark mixing in the Cabibbo--Kobayashi--Maskawa (CKM) matrix 
\cite{Cabibbo,KM} seems to imply).
There are certainly some subtleties here that may repay further study. 

In particular, there is not a single concept of
`weak doublet', 
but three separate concepts of `weak doublet', which can be presumably allocated as one for each of the generations, arising
from the three separate $Z_2$ quotients of $Q_8$. Each such quotient determines (and is determined by) a normal subgroup $Z_4$, that is the
kernel of the quotient map, so that there is also a correspondence between the three generations and the three $Z_4$ subgroups of $Q_8$,
generated by $i$, $j$ and $k$ respectively.

It may be helpful to take out only the positive real scalars from the group algebra, to leave a group
\begin{align}
Z_2 \times Z_2 \times Z_2 \times Z_2 \times SU(2)
\end{align}
which consists of the gauge group $SU(2)$ together with four quantum numbers.
It would be reasonable to assume that these are weak hypercharge and three components of weak isospin, but it is far from
obvious which basis of quantum numbers is most appropriate here. The finite group picks out a canonical basis, consisting of
one singlet, perhaps charge or weak hypercharge, and three doublets, perhaps weak isospin or neutral/colour charges.

If we are to attempt to predict some mixing angles from this group algebra, then we must use the geometry of $\HH$,
rather than the geometry of $\CC$ that we used in the previous two cases. If we include any leptons, we must include all three generations.
But we do not necessarily have to include all three generations of quarks. If we include charged leptons but not neutrinos, then we can choose one
combination of quarks (for example a proton) without the other half of the doublet (the neutron), and thereby span the $4$-space using only
three generations of electron, plus the proton. By omitting the neutrinos, and breaking the symmetry between proton and neutron,
we are effectively talking about the `right-handed' particles. We can then define a mass axis within the quaternions, just as we defined a mass axis in the
complex numbers previously. 

In the group algebra, of course, $\HH$ represents the left-handed particles, and is coordinatised with weak hypercharge and weak isospin.
The rest of the algebra is the tensor square of $\HH$
\begin{align}
\HH \otimes _\HH \HH = \RR + \RR + \RR + \RR,
\end{align}
in which the broken symmetry allows us to treat the three generations independently, each with its own mass parameter.
The quaternionic tensor product amounts, in principle, to nothing more than ordinary quaternion multiplication, so that the right hand side can also be
written naturally as a quaternion. However, we must be very careful to distinguish correctly between left-multiplication and
right-multiplication. If we are ever required to convert between left-multiplication and right-multiplication, then we must
at the same time take quaternion conjugates where necessary.

\section{Towards a quantisation of mass}
\label{massq}
\subsection{Quantum numbers}
At this stage we need some quantum numbers that generalise the $2$-space of quantum numbers generated by $Y_W$ and $T_3$. 
But if we are talking about `right-handed' particles, for which $T_3=0$ and $Y_W=2Q$, 
weak hypercharge and weak isospin are not useful, so we need a different set of quantum numbers altogether.
Note that we are not aiming to provide quantum numbers for all particles, but only for electrons in three generations, and the proton.

An alternative to using weak hypercharge and weak isospin
is to use the electric charge $Q=T_3+\frac12Y_W$ and neutral charge $N=T_3-\frac12Y_W$, as defined on the
`left-handed' electrons, and to transfer \emph{both} types of charge to the `right-handed' electrons.
This is not what is done in the Standard Model, but we shall see that it leads to a useful quantisation of mass,
once we split the neutral charge into three components, one
for each of the three generations. 
Instead of a single set of four quantum numbers, as in the Standard Model, we effectively have two sets: one in $\HH$,
and the other in $\HH\otimes_\HH\HH$. Logically, the former contains the weak hypercharge and weak isospin, while the
latter contains the electric and neutral charges.

The neutral charge $N$ takes the value $0$ on the `left-handed' electron, so let us assume it is also
$0$ on the `right-handed' electron.
Writing these charges in quaternionic notation, it makes sense to take the real part $1$
as the electric charge, and the three imaginary parts $i$, $j$ and $k$ as different types of neutral charge. 
Let us put the $0$ component of the neutral charge of the first-generation electron in the $k$ position, by analogy with the
third component of weak isospin.

At least one of the other components of neutral charge must be non-zero, and symmetry suggests they should be equal,
so let us scale this value to $1$. 
Rotating $i,j,k$ corresponds to a generation symmetry on leptons, which leads us to the following choice of quantum numbers
\begin{align}
e_R & =-1+i+j\cr
\mu_R & =-1+j+k\cr
\tau_R & =-1+i+k.
\end{align}
The proton must surely be symmetric in $i$, $j$ and $k$, which almost forces 
\begin{align}
p=1+i+j+k.
\end{align}

From these assumptions we can calculate the neutral combination
\begin{align}
e+\mu+\tau+3p = 5(i+j+k).
\end{align}
Now it was shown in \cite{remarks} that the empirical mass formula
\begin{align}
\label{emutau}
m(e)+m(\mu)+m(\tau)+3m(p) = 5m(n)
\end{align}
holds to well within current experimental uncertainty,
so that it makes sense to identify the neutron as
\begin{align}
n=i+j+k,
\end{align}
in order for the mass of these five particles to be obtained as a linear combination of the four quantum numbers.

Indeed, this proposed identification of the neutron confirms that our choice of quantum numbers is consistent.
Moreover, the lack of symmetry between protons and neutrons in the mass equation confirms that 
these quantum numbers apply to `right-handed' particles only.
Indeed, it may 
make more sense to apply the terms `left-handed' and `right-handed' to the quantum numbers,
rather than the particles, so that in this sense weak isospin is `left-handed', describing $\HH$, 
and the four charges are `right-handed', describing $\RR+\RR+\RR+\RR$.

Alternative coordinates can be obtained by multiplying by $1+i+j+k$, to give
\begin{align}
e & = -3+i-j-k\cr
\mu & = -3-i+j-k\cr
\tau & = -3-i-j+k\cr
n & =-3+i+j+k
\end{align}
so that the imaginary part then corresponds to the notation with $e,f,g$ in Section~\ref{toyelectrons}.
The real part now represents the (third component of) weak isospin, scaled by a factor of $6$. 
A further multiplication by $i+j+k$ gives
\begin{align}
e & = 1 -3i-j-5k\cr
\mu & = 1 -5i-3j-k\cr
\tau & = 1-i-5j-3k\cr
n & = -3 -3i-3j-3k
\end{align}
Then the weak isospin is measured in the direction $i+j+k$, and is the same for all four of these particles, as it should be. The real part is
a combination of charge and weak hypercharge, specifically $2Q-3Y_W$, so that these coordinates are compatible with those in
Section~\ref{Z4symmetry}.

\subsection{Calculating the mass vector}
In the earlier example
(Section~\ref{massplane}), we used the mass ratio of $\tau-\mu$ to $\tau-e$ as the single dimensionless parameter to define the single mixing angle,
apparently $\theta_{12}$ from the PMNS matrix.
Here we need three mass ratios to define three mixing angles. If we use the mass ratios of the three generations of electron and the proton, then
we would only expect to obtain two new (independent) mixing angles. However, it is also possible that the dependencies 
(if any) between mixing angles are
complicated four-dimensional geometrical relationships that are far from obvious, 
so that we may not notice them even if they are there.

We 
can define a mass vector 
\begin{align}
m & = m_0+m_1i+m_2j+m_3k
\end{align}
by substituting in the quantum numbers to obtain the equations
\begin{align}
m(n) & = m_1+m_2+m_3\cr
m(p) & = m_0+m_1+m_2+m_3\cr
m(e) & = -m_0+m_1+m_2\cr
m(\mu) & = -m_0+m_2+m_3
\end{align}
Taking these four independent masses from CODATA2018
\cite{CODATA2018},
in units of eV/$c^2$, to the nearest integer
\begin{align}
n & = 939565420\pm1\cr
p & = 938272088\cr
e & = 510999\cr
\mu & = 105658375\pm 3
\end{align}
we can solve the equations to obtain (up to an arbitrary scale factor) 
\begin{align}
m = -1293332 + 835200377i - 835982710j + 940347753k
\end{align}

\subsection{Mixing angles from geometry}
We now need to choose three relevant angles between this mass vector and elements of the regular configuration. As in the previous
examples, the angles are given by trigonometric formulae that are ratios of sums and differences of coordinates of the mass vector.
There are only three independent parameters, but a multitude of choices for which parameters we consider to be primary.

The three that 
appear to give the 
clearest translation to the standard model are
\begin{align}
(m_0-m_1-m_2)/m_0 & \approx .395103\cr
(m_1-m_2-m_3)/m_3 & \approx .777196874\cr
(m_1+m_2+m_3)/m_3 & \approx .999168039
\end{align}
The numerators come from the electron (first case) and the electron times $k$ (second case) and the neutron (third case).
The denominators are whichever of the real and $k$ parts actually occurs in each case.
The choice of $k$ amounts to a choice of generation of electron, and 
implies a choice of complex structure within the quaternions.
It is an obvious choice within the chosen coordinate system, though perhaps
not quite so obvious in the standard coordinate system of weak hypercharge and weak isospin.

We can also express these ratios in terms of the particle masses via
\begin{align}
m_0 & = m(p) - m(n)\cr
m_3 & = 2m(n) -m(p)-m(e)\cr
m_0 -m_1-m_2 & = - m(e)\cr
m_1+m_2+m_3 & = m(n)\cr
m_1-m_2-m_3 & = 3m(n)-2m(p)-2m(\mu)
\end{align}
so that
\begin{align}
\frac{m_0-m_1-m_2}{m_0} & = \frac{m(e)}{m(n)-m(p)}\cr
\frac{m_1-m_2-m_3}{m_3} & = \frac{3m(n)-2m(p)-2m(\mu)}{2m(n)-m(p)-m(e)}\cr
\frac{m_1+m_2+m_3}{m_3} & = \frac{m(n)}{2m(n)-m(p)-m(e)}
\end{align}
The mass ratios here are in most cases of \emph{differences} of particles, as they were in the previous examples. 
In the first case, both numerator and denominator have total charge $-1$, while in the other two cases they have total charge $0$.
The first is an obvious mass ratio to consider, the other two perhaps slightly less obvious

Of these three, the second is remarkably close to the value
\begin{align}
1/2+1/\sqrt{13} & \approx .777350098
\end{align}
that we have already seen, so probably does not give us anything new, just another manifestation of the electro-weak mixing angle.
It is hardly a surprise that the same angle should appear in the group algebra of $Q_8$ that appears in the subalgebra
already considered, that is the group
algebra of the subgroup $Z_4$.
The difference between the two values is around $.02\%$, while the given $9$-figure accuracy is only possible because we have used
the conjectured equation (\ref{emutau}) in an essential way.

If we do not use this equation, then the standard relative uncertainty in the tau mass of around $.01\%$ is large enough that 
we have to wonder whether the discrepancy
is statistically significant. In fact, re-calculating the mass ratio using the experimental $\tau$ mass of $1776.86\pm.12$ MeV$/c^2$ gives us a value
of $.777199\pm.000037$, which locates the discrepancy at around $4\sigma$.

An alternative 
interpretation is that  
the two angles are actually physically different, and their remarkable closeness is purely a coincidence.
This would mean that the value of the Weinberg angle used in calculating $\sin^2\theta_W$ for use with fermions is actually different
from the value used in $\cos\theta_W$ to determine the mass ratio of the W and Z bosons. Indeed, this is actually a \emph{prediction} of the
model presented here, since experimental uncertainty in the $e$, $\mu$, $p$ and $n$ masses is several orders of magnitude too small
to explain the discrepancy. For reference, the fermionic and bosonic values of the angle are
\begin{align}
\theta_W^f & \approx 28.154966^\circ\cr
\theta_W^b & \approx 28.165516^\circ
\end{align}

The other two mass ratios are probably best interpreted as cosines of angles: 
\begin{align}
.395103 & \approx \cos (66.7276^\circ)\cr 
.999168039 &\approx \cos (2.337325^\circ).
\end{align}
The latter relates to the neutron, therefore to quarks, hence to the CKM quark-mixing matrix, but fixes the first generation
(up and down quarks), so can only
plausibly be the mixing angle between the second and third generations, experimentally determined in the range 
\begin{align}
\theta_{23} = 2.38^\circ\pm.06^\circ.
\end{align}
The former can be written as
\begin{align}
(m(e)-m(\nu))/(m(n)-m(p)),
\end{align}
since the neutrino mass is too small to make a material difference to the equation, and can then be interpreted as
a relationship between the mass-difference-to-charge ratios for (first generation) leptons and baryons.
If it is a mixing angle, then it mixes the electroweak theory, that applies to both leptons and quarks, with quantum chromodynamics
(QCD), that applies only to quarks (and therefore baryons).
The only available mixing angle is therefore the CP-violating phase in the CKM matrix,
which is experimentally determined in the range
\begin{align}
\cos(68.8\pm 4.5)^\circ & \approx .362 \pm .075.
\end{align}
Again we have close agreement between prediction and 
experiment. 

\subsection{Taking stock}
At this stage we have plausible derivations of four of the nine mixing angles in the standard model. One of these uses
only electric charge, so gives an exact value based only on the fundamentals of the quark model. The other three use
measured values of (four) 
particle masses, and are therefore dependent to some extent on experiment, and in particular on the
calibration of mass. In all cases, the essential ingredient in the calculation is the geometry of the group algebra defined by
the finite symmetry groups. The symmetries are in all cases broken by either the charge (in one case) or the mass (in the other three cases).

In fact, we have \emph{two} plausible derivations of the electro-weak mixing angle, which give slightly different values, that are too close
together to be distinguishable by experiment. One (fermionic) value is obtained using charge and weak hypercharge alone, with a single generation
and without mass. The other (bosonic) 
value is obtained using mass and all three generations. It is possible, therefore, that distinguishing these two values
may extend the standard model in a useful way. 
Essentially this means breaking the connection between the fermionic usage of
$\sin^2\theta_W$ in the vertex factors of the Feynman calculus, and the bosonic usage of $\cos\theta_W$ as the mass ratio
of the W and Z bosons. This suggestion raises all sorts of questions about the nature of mass, and the Higgs mechanism for generating mass.

Perhaps more important than all of this, however, is the introduction of new quantum numbers for three `generations' of `neutral charge',
that provide a basis on which to quantise mass (in some instances)
as a four-dimensional concept, rather than the conventional scalar (one-dimensional) concept.
It might indeed be reasonable to rename these neutral charges as `gravitational' charges.

\section{The binary tetrahedral group}
\label{tetrahedron}
\subsection{Combining the two finite groups}
\label{overview}
The fact that the three generations of electrons appear in both the analysis of $Z_3$ and in the analysis of $Q_8$ suggests that we should 
combine these two symmetry groups to get a version of the weak $SU(2)$ in which the generation symmetry group $Z_3$ is explicit.
This means adjoining an automorphism $f$ of order 3 to the quaternion group, identifying $e=1$ and $g=f^{-1}$, and adding in the
new rules:
\begin{align}
if=fj,\quad jf=fk,\quad (\Rightarrow  kf=fi).
\end{align}
The resulting group is well-known and much studied. It is called the binary tetrahedral group, and has order 24. Every element is either an
element $q$ of the quaternion group, or of the form $qf$ or $qg$.

An explicit copy of the group within the quaternion algebra can be obtained by adjoining to the quaternion group $Q_8$
the quaternion
\begin{eqnarray}
w&:=& (-1+i+j+k)/2,
\end{eqnarray}
corresponding to $f$, and therefore also its quaternion conjugate
\begin{eqnarray}
v&:=& (-1-i-j-k)/2,
\end{eqnarray}
corresponding to $g$.
This group and its representation theory are described in some detail in \cite{perspective}, but we shall require 
yet more detail here.

The structure of the group algebra is well known (see \cite[p.404]{JamesLiebeck} and/or \cite{perspective,tetrions}): 
\begin{align}
\RR + \CC + M_3(\RR) + \HH + M_2(\CC),
\end{align}
where $M_n$ denotes the algebra of $n\times n$ matrices, with real or complex entries as specified.
In here we see the group algebra of $Z_3$ in the components $\RR+\CC$, while the group algebra of $Q_8$ consists of the
components $\RR+\HH$, plus the diagonal matrices in $M_3(\RR)$. The rest of $M_3(\RR)$ acts by `mixing' the three generations.
Indeed, the six off-diagonal terms can be split into three symmetric matrices and three anti-symmetric, 
the latter being generators for the compact subgroup $SO(3)$, in which we see precursors for
three generation-mixing parameters from the CKM and/or PMNS matrices.
Thus $M_3(\RR)$ encodes some of the properties of the strong force, although it does not actually contain the gauge group $SU(3)$.

Taking out four real scalars and one complex scalar from the real group algebra leaves us with the group
\begin{eqnarray}
U(1) \times SU(2) \times SL(2,\CC)\times SL(3,\RR),
\end{eqnarray}
which contains all the groups we need for the standard model of particle physics, apart from the fact that we have the
split real form $SL(3,\RR)$ rather than the compact real form $SU(3)$ of the Lie group of type $A_2$.

This difference may reflect the fact that we are looking at finite (generation) 
symmetries that can be observed in
experiments, rather than continuous
(colour) symmetries that are not observable. It may be possible to resolve this
issue by extending to the complex group algebra, but this process obscures a number of important features of
the real group algebra, so I will only do this if absolutely necessary. 
In any case, we must bear in mind that this difference must be resolved at some point,
since it may be the difference between a viable model and an unviable model.

\subsection{Irreducible representations}
\label{irreps}
Representation theory is usually presented first over the complex numbers, since this is the simplest theory, but we shall
require the extra subtlety of the representation theory over real numbers.
First note that the conjugation action of the group on itself divides the group into $7$
conjugacy classes, of sizes $1,1,6,4,4,4,4$, as listed in this table: 
\begin{eqnarray}
&\begin{array}{cccc}
\mbox{Size} & \mbox{Elements} & \mbox{Order}\cr\hline
1 & 1 = e& 1\cr
1 & -1 = i^2& 2\cr
6 & \pm i, \pm j, \pm k & 4\cr
4 & w, wi, wj, wk & 3\cr
4 & v, -vi,-vj, -vk & 3\cr
4 &- w, -wi, -wj, -wk & 6\cr
4 & -v, vi,vj, vk & 6\cr\hline
\end{array}&
\end{eqnarray}

There are therefore $7$ complex irreducible representations, whose characters are listed,
following the usual mathematical convention, as rows in the following table. The top row gives a representative of the conjugacy class, 
the entries are the traces of the representing matrices, 
and $\omega$ and $\bar\omega$ are the complex cube roots of unity.
\begin{eqnarray}
\begin{array}{ccccccc}
1&-1&i&w&-w&v&-v\cr\hline
1&1&1&1&1&1&1\cr
1&1&1&\omega&\omega&\bar\omega&\bar\omega\cr
1&1&1&\bar\omega&\bar\omega&\omega&\omega\cr
3&3&-1&0&0&0&0\cr
2&-2&0&-1&1&-1&1\cr
2&-2&0&-\omega&\omega&-\bar\omega&\bar\omega\cr
2&-2&0&-\bar\omega&\bar\omega&-\omega&\omega\cr\hline
\end{array}
\end{eqnarray}
The three $2$-dimensional representations will eventually play the role of Weyl spinors, but there are three of them, 
not just the familiar left-handed and right-handed Weyl spinors. 
 Of course, the standard model for electro-weak interactions really has three Weyl spinors, two left-handed and one right-handed, so the question is whether it is possible, and if so, how, to match up the three spinors in the two models.

The structure of the complex group algebra 
can be read off from the dimensions of the representations, and is
\begin{align}
3\CC + M_3(\CC) + 3M_2(\CC),
\end{align}
consisting of one copy of $n\times n$ complex matrices for each complex $n$-dimensional representation.
In addition to the problem of mixing and matching the three copies of $SU(2)$ or $SL(2,\CC)$ inside the three copies of $M_2(\CC)$,
 there is the problem of mixing and matching the three copies of $U(1)$ inside the three copies of $\CC$. There is too much choice at this stage,
 and the problem becomes rather more tractable if we restrict to the real group algebra. 
 
Moreover, as we have seen in the toy models considered above, the real group algebra distinguishes real, complex and quaternionic matrix algebras,
which has already proved useful for deriving particular parts of the Standard Model.
 Of course, this restriction means that we lose the group $SU(3)$
from $M_3(\CC)$, but we have $SL_3(\RR)$ instead. This different real form may or may not be a satisfactory replacement for the 
gauge group of the strong force. Indeed, there is also a spare copy of $U(1)$ inside the real group algebra, so that we can use $U(1)\times SL(3,\RR)$
in place of $SU(3)$ if this helps. The basic idea would be to create the required $3$-dimensional complex representation
as the tensor product of a $1$-dimensional complex representation of $U(1)$ and a $3$-dimensional real representation
of $SL(3,\RR)$.

The irreducible real representations of the binary tetrahedral group
can be described by taking the sum of a complex representation with its complex conjugate, and at the same time
taking the union of each conjugacy class with its inverse class. This reduces the table to five rows and five columns, as follows:
\begin{align}
\begin{array}{ccccc}
1&-1&i&w&-w\cr\hline
1&1&1&1&1\cr
2&2&2&-1&-1\cr
3&3&-1&0&0\cr
4&-4&0&-2&2\cr
4&-4&0&1&-1\cr\hline
\end{array}
\end{align} 

The first of the two $4$-dimensional representations is quaternionic, so has a continuous $3$-parameter family of possible complexifications.
The second of the two $4$-dimensional representations, on the other hand, has a discrete set of exactly two possible complexifications. This interplay between the discrete and the continuous 
may possibly play an important role in the structure of the Dirac spinor and the Dirac algebra \cite{Dirac}.
In order to distinguish these two representations, let us write the first one $4_H$, to denote a Hamiltonian quaternionic representation, 
and the second one $4_C$, to denote a classical complex representation. Then the group algebra, 
as a representation of the finite group $G$, has the structure
\begin{align}
1+2+3+3+3+4_H+4_C+4_C.
\end{align}

\subsection{Tensor products}
\label{tensors}
So far, I have associated the various Lie groups (gauge groups and spin groups) with a single representation of the finite group.
The finite group links them together in various ways, described in part by the decompositions of tensor products of representations.
These decompositions can be calculated easily from the character table, and are as follows, with $+$ signs omitted to save space:
\begin{align}
\begin{array}{c|cccc}
1&2 & 3 &4_H & 4_C\cr\hline
2 & 13 & 123 & 4_C4_C & 4_H4_C\cr
3 & 123
& 1233 & 4_H4_C4_C & 4_H4_C4_C\cr
4_H & 4_C4_C&4_H4_C4_C & 11113333 & 223333\cr
4_C & 4_H4_C&4_H4_C4_C & 223333
& 1123333 
\end{array}
\end{align}

There is clearly plenty of structure in here that might relate to the structure of the standard model. For example the fact that $(1+2)\otimes 4_C = 4_H+4_C+4_C$ means that the whole of the fermionic part of the algebra can be obtained from a real plus a complex `version' of the spin representation $4_C$.

 If one distinguishes the two copies of the real numbers by labelling one of them with $\gamma_5$, then there is 
 some prospect of being able to match up the finite model with the standard model. Similarly, all of the $4\otimes 4$ representations look like different real slices of the complex Clifford algebra in the standard model, which for a finite (therefore compact) group must split as $1+(1+3)+(3+3)+(1+3)+1$. Again we see a mixing of the top and bottom degrees of the Clifford algebra ($1$ and $i\gamma_5$) into a real $2$-space, and a similar mixing in the odd part of the Clifford algebra. All this suggests that a careful distinction between $4_H$ and $4_C$ may be able to throw some light on the mixing of quantum electrodynamics and the weak interaction in the standard model. 

A similar picture emerges from the bosonic part of the algebra, in which one obtains the whole algebra from the tensor product $(1+3)\otimes 3 = 1+2+3+3+3$. For the purpose of matching with the standard model, one might wish to go further and observe the isomorphism
\begin{align}
(1+3)\otimes (1+3)&\cong1+1+2+3+3+3+3\cr
&\cong 4_C\otimes 4_C. 
\end{align} 
However one tries to interpret this equivalence, it seems to imply a mixing between the strong force, with a gauge group acting on $1+3$, and the
electroweak forces, with a gauge group acting on $4_C$.

Indeed, there are a number of other suggestive isomorphisms between different tensor product representations. For example, for any
representation $R$ that does not involve the trivial representation $1$, we have
\begin{align}
(1+2)\otimes R &\cong  3\otimes R,
\end{align}
which gives ample scope for mixing a broken $1+2$ symmetry with an unbroken $3$ symmetry.
Another example is that
\begin{align}
3\otimes 4_H &\cong
3\otimes 4_C,
\end{align}
which gives plenty of scope for mixing a spinor of type $4_H$ with a spinor of type $4_C$.

\subsection{Explicit matrices}
\label{matrices}
For the purposes of explicit calculation, it will be useful to have explicit matrix copies of all the irreducible representations. It is sufficient to specify matrices representing the generators $i$ and $w$. In the $1$-dimensional representation, they both act trivially. In the $2$-dimensional representation, $w$ is a rotation of order $3$, and $i$ acts trivially, so we may take the matrices
\begin{eqnarray}
i\mapsto \begin{pmatrix}1&0\cr 0&1\end{pmatrix}, && w\mapsto \frac12\begin{pmatrix}-1& \sqrt 3\cr -\sqrt 3 & -1\end{pmatrix}
\end{eqnarray}
The $3$-dimensional representation is the representation as symmetries of a regular tetrahedron, which can be embedded as alternate vertices of a cube, so that the matrices can be taken as
\begin{eqnarray}
i\mapsto\begin{pmatrix}1&0&0\cr 0&-1&0\cr 0&0&-1\end{pmatrix},&&
w\mapsto\begin{pmatrix}0&1&0\cr 0&0&1\cr 1&0&0\end{pmatrix}.
\end{eqnarray}

The representation $4_H$  is the representation by right-multiplication on the quaternions, so that the matrices can be taken as
\begin{eqnarray}
i\mapsto \begin{pmatrix}0&1&0&0\cr -1&0&0&0\cr 0&0&0&-1\cr 0&0&1&0\end{pmatrix},&&
w\mapsto\frac12\begin{pmatrix}-1&1&1&1\cr -1&-1&-1&1\cr -1&1&-1&-1\cr -1&-1&1&-1\end{pmatrix}
\end{eqnarray}
If this representation corresponds to a spinor in the standard model,
then the quaternionic symmetry is broken, and a particular complex basis must be chosen. The following matrices give an example,
which may or may not be similar to what is done in the standard model:
\begin{eqnarray}
i\mapsto \begin{pmatrix} i & 0\cr 0 & -i\end{pmatrix},&& w\mapsto \frac12\begin{pmatrix}-1+i&1+i\cr -1+i & -1-i\end{pmatrix}
\end{eqnarray}

Finally, the representation $4_C$ can be written as 
\begin{eqnarray}
i\mapsto \begin{pmatrix}0&1&0&0\cr -1&0&0&0\cr 0&0&0&-1\cr 0&0&1&0\end{pmatrix},&&
w\mapsto \begin{pmatrix}1&0&0&0\cr 0&0&1&0\cr 0&0&0&1\cr 0&1&0&0\end{pmatrix}
\end{eqnarray}
It is worth pausing for a moment to compare this representation with $1+3$: the action of $w$ is the same in both cases,
but the actions of $i,j,k$ are not.

\subsection{Projections}
\label{projections}
The theory of finite group representations contains a canonical set of
projections from the group algebra onto the
various matrix subalgebras. 
Mathematically, a projection $p$ is an \emph{idempotent}, which just means that $p^2=p$.
Two idempotents $p$ and $q$ are \emph{orthogonal} if $pq=0$, so that 
\begin{align}
(p+q)^2 & =p^2+q^2\cr
& =p+q
\end{align}
 and
$p+q$ is also an idempotent. Any idempotent that cannot be split is a sum of two othogonal idempotents is called \emph{primitive}.
The primitive idempotents in the group algebra are equivalent to the identity elements in the matrix algebras.

First there is the basic division into bosons and fermions, obtained via
projection with the idempotents $(e\pm i^2)/2$. 
Within the fermions, there are two projections obtained from the
sum of all the elements of order $3$, that is
\begin{align}
s&:= f+ if+ jf+ kf + g-ig-jg- kg\cr
& = (1+i+j+k)f + (1-i-j-k)g.
\end{align}
These projections are defined by the elements $(2e-s)/6$ and $(4e+s)/6$ of the group algebra. Notice that the coefficients
of the identity element here are $1/3$ and $2/3$, which suggests some connection with the charges on the quarks.

Within the bosons, 
there are three projections, so that the full list is as follows:
\begin{eqnarray}
&\begin{array}{cc}
\mbox{Subalgebra} & \mbox{Idempotent}\cr\hline
\RR & (e+i^2)(e+i)(e+j)(e+f+g)/24\cr
\CC & (e+i^2)(e+i)(e+j)(2e-f-g)/24\cr
M_3(\RR) & (e+i^2)(3e-i-j-k)/8\cr\hline
\HH & (e-i^2)(2e-s)/12\cr
M_2(\CC) & (e-i^2)(4e+s)/12 \cr\hline
\end{array}&
\end{eqnarray}

If we omit the subalgebra $M_3(\RR)$ from this description 
we obtain the $15$-dimensional subalgebra
\begin{align}
&\RR+\CC+\HH + M_2(\CC)
\end{align}
 in which there is a close parallel between the splitting of the `scalars'
 $\RR+\CC$ and the splitting  of the `spinors'
 \begin{align} 
 \HH + M_2(\CC) = \HH \otimes (\RR+\CC)
 \end{align}
 using the different splittings of the identity element as
 \begin{align}
  e &= (e+f+g)/3 + (2e-f-g)/3,\cr  
  e&= (2e-s)/6 + (4e+s)/6.
 \end{align}

These two  
pairs of projections appear to be the finite analogues of the 
pair of projections $(1\pm \gamma_5)/2$
that split the Dirac spinor into left-handed and right-handed Weyl spinors in the standard model,
and that distinguish the weak force from electromagnetism. 
However, the fact that in the finite model the projections on the scalars are different from the projections on the spinors
allows the finite model to be more subtle than the standard model, and perhaps incorporate a generation structure
into this distinction.  It may be possible, for example, to add in the scalars from $M_3(\RR)$, to obtain a $16$-dimensional
algebra that fulfils the function of the Dirac algebra, but with a more subtle structure derived from the action of the finite group.

It is worth remarking here that because these projections involve dividing by $2$ and $3$, they can be implemented in the real group algebra
$\RR G$ or the rational group algebra $\QQ G$, but not in the integral group ring $\ZZ G$. Indeed, the structure of
$\ZZ G$ is much more subtle than that of $\RR G$. Ultimately, to implement the finite model in full we will need to grapple
with this structure in detail. For the purposes of the present paper, however, it is enough consider only the real group algebra. 

Just to give a flavour of the implications of using the integral group ring, I'll describe a toy model of the proton using the integral group
ring of the group $Z_3$ of order $3$. If we take $e,f,g$ as the elements of $Z_3$, then the idempotents in $\RR Z_3$ are
$(e+f+g)/3$ and $(2e-f-g)/3$ as above. The former projects onto a real $1$-space containing a down quark, and the latter projects onto a real
$2$-space containing two up quarks. These are perfectly good real representations of $Z_3$, and give a perfectly good description of the
internal structure of a proton.

But they are not integral representations of $Z_3$. The integral 
group algebra $\ZZ Z_3$ does not support any projections, and is
indecomposable. This I interpret as saying that the proton itself cannot be decomposed into smaller particles, so that, in the real
universe, protons never decay. Thus the quarks are `real' particles, but they are not `whole' (integer) particles.

It is also possible to use the action of $Z_3$ by conjugation to describe the corresponding bosons. These bosons consist of one `left-handed'
and one `right-handed' quark, that correspond to `quark' and `anti-quark' in the standard model. So I get three pions this way, consisting of two 
charged pions, $\pi^+=u\bar d$ and $\pi^-=d\bar u$, and a neutral pion $\pi^0=u\bar u$. While this description of the neutral pion is not the same
as in the standard model, it is the only possibility in a discrete model, in which quantum superposition cannot be implemented as an intrinsic property
of an elementary particle, but only as a property of the experiment or the environment.

\section{Electroweak unification}
\label{EWunif}

\subsection{The Lorentz group} 
The final component of the algebra is $M_2(\CC)$, which consists of a real scalar, a second copy of $U(1)$, and a group $SL(2,\CC)$.
The last would normally be interpreted as the Lorentz group, acting by left-multiplication on a Dirac spinor, so that a Dirac spinor 
effectively consists of two Weyl spinors corresponding to the two columns of an element of $M_2(\CC)$. There are some technicalities
concerning the additional copy of $U(1)$, and complex conjugation, that appear in the standard model, 
which will require detailed investigation at some point. In particular, the two columns of $M_2(\CC)$ naturally have the same (Weyl)
handedness, rather than opposite handedness, so that the translation to the standard model is not necessarily straightforward.

The important point to note here is that the Lorentz group \emph{emerges} from this discrete form of electro-weak unification, 
using the groups $Z_3$ and $Q_8$ in place of $U(1)$ and $SU(2)$, and does not have to be added
as a `background'.
Algebraically, we take the linearised gauge groups
\begin{align}
\RR U(1) & = \CC\cr
\RR SU(2) & = \HH
\end{align}
and construct the tensor product of the resulting algebras
\begin{align}
\CC \otimes_\RR \HH = M_2(\CC).
\end{align}

Of course, the Lorentz group does \emph{not} arise automatically if only one of the two forces is considered. 
Hence it is necessary
to construct the Dirac algebra separately from the gauge group in the standard theory of quantum electrodynamics (QED). 
Subtle properties of spacetime cannot therefore be probed by either electromagnetic or weak interactions on their own, but only
by the combination of the two. 

This underlines the great importance of the famous Wu experiment \cite{Wu} of 1957, the essence of which was to
combine an electromagnetic property (spin) with a weak property (beta decay), and to measure the results with respect to the locally defined
directions of spacetime. The result of the experiment was to confirm a correlation between the three representations. The electromagnetic
symmetries distinguished between spin up and spin down, the weak symmetries distinguished between the emitted electron and antineutrino,
and the spacetime symmetries distinguished left-handed from right-handed frames of reference. The experiment proved that these three
dichotomies are not independent: 
there is a correlation between the direction of spin and the direction of momentum of the ejected electron.

\subsection{Chirality}
The group algebra hypothesis does much the same thing, theoretically instead of experimentally, via an embedding of $Z_4$ in $Q_8$. 
There are three such embeddings, one for each generation, so to analyse the Wu experiment we must choose the
first generation copy. The quotient group $Q_8/Z_4$ then acts on the weak doublets, and effectively swaps the electron with the
antineutrino. In terms of Fig.~\ref{Wfig}, the action is a reflection on the $2$-space, and can therefore be implemented as
complex conjugation on the complex representation of $Z_4$. This complex conjugation happens naturally within the
subalgebra $\HH$ of the group algebra of $Q_8$, and also in the group algebra of the binary tetrahedral group. 

In the latter case, $Q_8$ also acts on the left-handed and right-handed Weyl spinors, in such a way that the eigenvectors of the
subgroup $Z_4$ are the spin up and spin down states, and $Q_8/Z_4$ swaps these two states.
In other words, the existence of the finite symmetry group implies that reversing the spin cannot be achieved without at the same time
swapping the electron with the antineutrino. Thus the group algebra explains the chirality of beta decay, as detected by the
Wu experiment. However, it does this by distinguishing spin up from spin down, and \emph{not} by distinguishing left-handed
spinors from right-handed spinors.

The latter distinction is made by complex conjugation on the subalgebra $M_2(\CC)$, which is equivalent to complex conjugation on
the subalgebra $\CC$, and cannot be achieved within the group algebra itself. It can however be achieved by an outer automorphism of the
finite group. There is therefore a big difference between complex conjugation on the representation of $Z_4$, and complex conjugation
on the representation of $Z_3$. The former is a physical symmetry, that can be realised by a rotation in space in order to reverse
the direction of the magnetic field, while the latter is an unphysical symmetry, that involves reversing the direction of time.

In the standard model, however, both are represented by complex conjugation, which makes it 
easy to confuse the two if one is not careful.
The group algebra approach clarifies this distinction, and makes clear that the \emph{chirality} of beta decay,
that is a property of $Z_4$ and the weak interaction,
has nothing to do with the (relativistic) distinction between left-handed and right-handed spins or \emph{helicities},
that is a property of $Z_3$ and electrodynamics. The discreteness, in other words, makes a very clear distinction between
two distinct copies of $U(1)$ that might otherwise be confused with each other.

\subsection{The remaining mixing angles}
So far I have identified one of the standard model mixing angles in $\CC$, one in $\HH$, and three in $M_3(\RR)$.
This includes two different versions of the weak mixing angle, one bosonic and one fermionic. In particular, I have identified all four
of the mixing angles that are available in the compact part $U(1) \times SO(3)$ of the bosonic half of the algebra.
The fermionic half has compact part $SU(2)\times U(2)$, so could in principle give us seven mixing angles in total.
However, the scalar part of $U(2)$ is not acted on by the finite group in the adjoint representation, so that it does not participate in
interactions, and therefore there is no angle there that can be physically measured. Thus there should be only six fermionic mixing angles,
one of which has been already identified, leaving exactly five still to identify.
This is exactly the number of remaining mixing angles in the Standard Model.
Explicitly, these are the angles $\theta_{12}$ (the Cabibbo angle) and $\theta_{13}$ from the CKM (quark-mixing) matrix,
together with $\theta_{13}$, $\theta_{23}$ and $\delta_{CP}$ from the PMNS (lepton-mixing) matrix. 

The mechanism for converting bosonic to fermionic mixing angles is to tensor with  the representation $\HH$, which has the effect
of extending from one generation to three. Hence from $\CC\otimes\HH = M_2(\CC)$ we obtain the three remaining
angles in the PMNS matrix, which puts the other two quark-mixing angles effectively in $\HH$.
This allocation implies that in order to complete the CKM matrix we need to use the first and second components of weak isospin,
not just the (arbitrary)  third component as in the Standard Model.
Similarly, the rest of the PMNS matrix, that describes neutrino oscillations, requires us to consider all three components of spin,
rather than the (arbitrary) $z$ component as in the Standard Model.

Assuming that these five mixing angles can indeed also be found by similar geometrical arguments in the components of the group algebra,
it would appear that this group algebra is actually big enough to contain the whole of the standard model, despite the absence of an
explicit copy of $SU(3)$. If so, then the standard model $SU(3)$ must be made by a process of complexifying $SL_3(\RR)$, perhaps
using the second copy of $U(1)$ that exists in the group algebra. 
Such suggestions have been made before \cite{remarks,MDW}, in response to the same problem that $SU(3)$ itself does not arise naturally
from the proposed mathematical model. The question is quite subtle, and comes down to the question of whether it is really necessary to
complexify the Lie group, or whether it is enough to complexify the Lie algebra, or even to complexify only the representation spaces.

For example, the solution adopted in \cite{MDW} is to combine operators $A$ in $U(1)$ and $B$ in $SL(3,\RR)$ acting on a spinor $\Psi$
via the operation
\begin{align}
[A,[B,\Psi]]
\end{align}
rather than 
\begin{align}
[[A,B],\Psi].
\end{align}
This enables the construction of the operators of $SU(3)$ that are required in the Standard Model, without having an explicit copy 
of the Lie algebra of $SU(3)$ in the underlying real form of the Lie algebra of $E_8$.
A similar procedure can be adopted in the group algebra model, in which $U(1)$ is derived from $Z_3$ acting on a complex
$1$-space, and $SL(3,\RR)$ is derived from $A_4$ acting on a real 3-space. The tensor product of these two representations
is then a complex $3$-space on which generators for $SU(3)$ can be constructed.

\subsection{Real parameters}
We have seen that the $11$ compact dimensions give rise to ten mixing angles, of which nine are in the Standard Model and one is a new
`weak' mixing angle of about $.01055^\circ$, or $38.0''$, being the difference between the fermionic and bosonic
versions of the Weinberg angle. In the same way, the $13$ non-compact dimensions give rise to a further $12$
real parameters, which we would like to identify specifically. Four of them are ``pseudoscalars'', that is scalars associated with 
individual forces, so are probably best interpreted as coupling constants. The non-scalars are the five boosts in $SL(3,\RR)$
and the three boosts in $SL(2,\CC)$. The former are related to the strong force, and the latter to electromagnetism, so are
essentially the masses of the quarks and electrons respectively.

The reason for only having five quarks masses rather than six may be that the top quark does not hadronise,
which means that it does not participate in the strong force that glues hadrons together. We might, however, 
prefer use the remaining
boost in $GL(3,\RR)$ for the top quark mass, instead of the strong coupling constant.
If so, we might also want to use the remaining boost in $GL(2,\CC)$ for a fourth electromagnetic mass,
in which case the proton mass is the obvious candidate.

The calculations above imply that, at least in some cases, the mixing angles can be calculated from the masses and coupling constants.
Indeed, the mathematical structure of the algebra implies that this can be done in all cases. The exact form of the calculations is not
clear in all cases, but in principle this can all be worked out. This 
means that the total number of independent parameters is roughly halved compared to the Standard Model. But it does not give us much of a clue
as to where these remaining parameters come from.

\section{Implementing the standard model}
\label{SM}
\subsection{Resolving conflicts}
At this stage it is clear that there are some fundamental conflicts between the group algebra model presented here,
and the Standard Model as usually understood. These may be serious enough as to be deal-breaking, but the potential benefits of the
group algebra model in terms of providing fundamental explanations of unexplained structures and parameters are such that
a serious effort should first be made to resolve these conflicts. 

The most obvious problem is that the gauge group $SU(3)$ of the
strong force does not appear explicitly in the group algebra. The question is whether it can be effectively constructed and used
to replicate QCD within the group algebra context. Since the fundamental particles arise in the Standard Model from tensor products
of representations of various groups, the basic question is whether the $3$-dimensional complex representation of $SU(3)$ is itself
fundamental, or whether it can be broken down into a tensor product of more fundamental objects. For example, in \cite{MDW} it is
taken to be a tensor product of a one-dimensional complex representation of $U(1)$ with a 3-dimensional real representation
of $SL(3,\RR)$. Here we may wish to restrict further to the compact group $SO(3)$, or even to the finite group.

A second problem concerns the Dirac spinor, which again does not occur explicitly as a sum of left-handed and right-handed
Weyl spinors. Instead we have a pair of Weyl spinors of the same handedness in $M_2(\CC)$, together with a quaternionic (Majorana)
spinor with no canonical complex structure at all. The question is whether a choice of complex structure within the quaternions
is sufficient to support a version of the Dirac equation that is equivalent to the standard Dirac equation for a single generation.

\subsection{Emergence of the strong force}
We can, for example, tensor the linearised gauge group of the weak force with itself. At the level of algebras this just gives
\begin{align}
\HH \otimes _\HH \HH = \HH
\end{align}
so does not appear to give us anything new.
But at the level of representations of the finite group we get only a real representation on the right hand side, splitting into the
one-dimensional antisymmetric square (real part) and the three-dimensional symmetric square (imaginary part)
\begin{align}
\HH \otimes_\HH \HH = \RR + \RR^3.
\end{align} 
Tensoring with $\CC$ as the linearised gauge group $U(1)$ 
then gives us
\begin{align}
\CC \otimes_\RR(\RR + \RR^3) = \CC + \CC^3,
\end{align}
which allows us to use $\CC^3$ as a representation for the colours of QCD.  

It is then possible to impose a unitary structure on $\CC^3$ and impose a symmetry group $SU(3)$ acting on it,
but this is not a natural construction in the group algebra, and has some rather serious consequences, which are well-known.
One is that the action of $SU(3)$ on rows does not commute with the action of the rotation group $SO(3)$ on columns.
In other words, $SU(3)$ causes the fundamental structure of space to break down. In particular, it breaks the symmetry of space
by picking a preferred direction. This choice of direction within a $3$-space appears to have a multitude of consequences,
depending on the interpretation of the $3$-space. It breaks the symmetry of the weak interaction; it separates
the gravitational field from the electromagnetic field; and it specifies the direction of gravity. 

In the group algebra model, we allow the whole of $GL_3(\RR)$ to act on columns, but restrict to the finite group acting on rows.
The finite group then commutes with $GL_3(\RR)$, and is in effect a finite version of $SO(3)$ acting on rows. This group $SO(3)$
can be identified with the unbroken version of the weak $SU(2)$, and with the subgroup $SO(3)$ of $SU(3)$, and hence
provides a method of unifying the strong and weak nuclear forces with a common gauge group $SO(3)$.

By associativity of tensor products, this representation can also be written in terms of a tensor product of
spinors with isospinors. However, we have to be careful  about whether we are tensoring over real numbers, complex numbers, or
quaternions. Since the spinors are complex representations, we must tensor over complex numbers, which means that,
in order to get the dimensions to work out correctly, we must (a) use a Weyl spinor instead of the whole Dirac spinor, and (b) break
the quaternionic symmetry of the isospinor to write it also as a complex two-dimensional representation:
\begin{align}
\HH \otimes_\CC \CC^2 = \CC + \CC^3.
\end{align}

Both (a) and (b) are conspicuous forms of symmetry-breaking in the Standard Model.
The former is responsible for the projections onto the left-handed and right-handed Weyl spinors, while the latter
implies a choice of a particular generation of electron. The group algebra model enables us to avoid both of
these problems, by re-associating the tensor product into the form
\begin{align}
\CC \otimes_\RR (\HH\otimes_\HH \HH) = \CC + \CC^3,
\end{align}
and making essential use of the quaternionic (rather than complex) tensor product in order to unify the three generations.

\subsection{The strong force from first principles}
The smallest group that has a %real 
$3$-dimensional representation is the tetrahedral group of order $12$. This group can be thought of 
in various ways, including (a) the rotation symmetry group of a regular tetrahedron, (b) the group of even permutations of four letters, 
(c) the quotient of the binary tetrahedral group by $Z_2$, obtained by ignoring the signs, and (d) the group $SL(2,Z_3)$ of all $2\times 2$ matrices of
determinant $1$ over the field $Z_3$ of three elements. 

Its group algebra is half of the group algebra of
the binary tetrahedral group, and has structure
\begin{align}
\RR + \CC + M_3(\RR).
\end{align}
Its representations therefore split into
\begin{enumerate}
\item the trivial representation, 
\item a complex $1$-space, and 
\item a real $3$-space. 
\end{enumerate}
The tensor product of (2) with (3)
is then a complex $3$-dimensional representation on which both $\CC$ and $M_3(\RR)$ act. It is therefore possible to construct the full
algebra of $3\times 3$ complex matrices acting on this complex $3$-space, and hence reproduce all of the necessary formalism of
quantum chromodynamics (QCD). However, the part of the unitary group $SU(3)$ that is reached by the finite symmetry group, that
describes physical interactions of particles, is only $SO(3)$, so that the rest of $SU(3)$ is in principle unobservable.
Hence 
replacing $SU(3)$ with $SL(3,\RR)$ cannot, even in principle, affect any Standard Model predictions of
experimental results.

\subsection{The Dirac gamma matrices}
The first and most fundamental requirement for relating the group algebra model to the Standard Model 
is to implement the Dirac equation. The most obvious, though not necessarily the best, 
place to do this is surely the matrix
algebra $M_2(\CC)$. Thus we must use these $2\times 2$ matrices to implement the Dirac spinor, divided into the left-hand
column and the right-hand column, that should be related in some way to the left-handed and right-handed Weyl spinors in the standard model.

Actions on the Dirac spinors can then be obtained from both left-multiplications and right-multiplications, which together generate a $16$-dimensional
complex algebra that we must identify with the complex Clifford algebra in the standard model. For this purpose, we need a choice of
Dirac gamma matrices. There are many possible choices, but since each of the four Dirac matrices swaps the left-handed and right-handed parts of the spinor, the following looks to be a good choice. Each Dirac matrix is written as a pair
of a left-multiplication and a right-multiplication by Pauli matrices.
\begin{align}
\gamma_0 =(1,\sigma_1)&= \left( \begin{pmatrix}1&0\cr 0&1\end{pmatrix}_L, \begin{pmatrix}0&1\cr1&0\end{pmatrix}_R\right),\cr
\gamma_1 =(\sigma_1,i\sigma_2)&= \left( \begin{pmatrix}0&1\cr 1&0\end{pmatrix}_L, \begin{pmatrix}0&1\cr-1&0\end{pmatrix}_R\right),\cr
\gamma_2 =(\sigma_2,i\sigma_2)&= \left( \begin{pmatrix}0&-i\cr i&0\end{pmatrix}_L, \begin{pmatrix}0&1\cr-1&0\end{pmatrix}_R\right),\cr
\gamma_3=(\sigma_3,i\sigma_2)&= \left( \begin{pmatrix}1&0\cr 0&-1\end{pmatrix}_L, \begin{pmatrix}0&1\cr-1&0\end{pmatrix}_R\right).
\end{align}

With these definitions, we can write down the Dirac equation in the usual way. The group algebra model, however, contains more structure than
the standard model, and in particular contains an important distinction between the Lie group acting on the left, and the finite group acting on the right.
Separating out the left-multiplications and the right-multiplications in the above, we have the following right-multiplications:
\begin{align}
i\gamma_0&=(1,i\sigma_1),\cr
i\gamma_1\gamma_2\gamma_3&=(1,i\sigma_2),
\end{align}
which together generate the quaternion group $Q_8$.
Hence this Dirac equation can be used to study four particles, as in the toy model discussed earlier.
The same equation applies to a variety of different particles, but only to four at a time. 
It is possible to use
this equation to study one generation of elementary fermions, but not to study all three simultaneously.
In particular, the standard model has to implement the generation structure outside the Dirac algebra.

On the left, we have 
\begin{align}
i&= (i,1)\cr
\gamma_1\gamma_2&=(i\sigma_3,1)\cr
\gamma_2\gamma_3&=(i\sigma_1,1),
\end{align}
which on exponentiation generate the Lie group $U(2)$, that is, a group isomorphic to the electro-weak gauge group in the standard model. 
Alternatively, we can obtain $SL(2,\CC)$ from $i\gamma_1\gamma_2$ and $i\gamma_2\gamma_3$.
Note, however, that this copy of $SL(2,\CC)$ is distinct from Dirac's relativistic spin group generated by
$\gamma_0\gamma_1$, $\gamma_1\gamma_2$ and $\gamma_2\gamma_3$. 
In particular, these two copies of $SL(2,\CC)$
have different physical meanings. 

\subsection{The Dirac equation}
As described above,
 every matrix subalgebra 
 of the group algebra is obtained by a canonical projection by a well-defined idempotent \cite{tetrions}, 
 so contains a canonical copy of the group.
 The rows can be scaled to any size, provided the columns are then scaled by the inverse. This allows any part of the theory to be
 (re)normalized to any desired scale. The idempotents themselves define natural units for the relevant physical concepts.
 For example, the three `bosonic' components $\RR+\CC+M_3(\RR)$ are connected by two fundamental constants that relate the three
 identity matrices to each other. These are the speed of light, and Planck's constant.
 Every differential or integral equation in fundamental physics arises from decomposing a tensor product of representations,
 via a well-defined process of quantisation, or de-quantisation (decoherence).
 
Canonical quantisation is usually expressed as converting a macroscopic momentum into a quantum of the dual of space via the operators
\begin{align}
i\hbar\partial_x, i\hbar\partial_y, i\hbar\partial_z.
\end{align} 
The factor of $i\hbar$ here is inserted into what is fundamentally a \emph{real}, not imaginary, operation, in order
to convert between the components $\RR$ and $\CC$. Similarly, a factor of $c$ is inserted into the energy term $\partial_{ct}$ in order
to link the components $\RR$ and $M_3(\RR)$. With these prerequisites, the Dirac equation arises from the fact that
the full `fermionic' part of the algebra has $6$ different decompositions as a tensor product of either $\RR+\CC$ or $\RR^3$
with either $\HH$ or a left-handed or right-handed Weyl spinor $\CC^2_L$ or $\CC^2_R$.

The simplest version comes from
\begin{align}
\RR^3\otimes \HH & = (\RR+\CC)\otimes \HH\cr
& = \HH + M_2(\CC)
\end{align}
and moving the energy term to the left-handed side, to get a projection
\begin{align}
(\RR+\RR^3) \otimes \HH & \rightarrow M_2(\CC).
\end{align} 
To convert this into the standard Dirac equation one has to mix $\HH$ with the left-handed part of $M_2(\CC)$ via a projection
with $1-\gamma_5$, in order to identify charge with a combination of weak hypercharge and weak isospin, 
but it is simpler just to use the above equation, which contains the same physical information. 
Keeping the weak isospin separate from the
left-handed Weyl spinor also allows us to make full use of the mass plane, and all three components of weak isospin, in order
to describe all the fundamental particles of matter at once, rather than one at a time.

\subsection{Electro-weak mixing}
In the standard model, electro-weak mixing \cite{Weinberg} is described as a mixing of the gauge groups $U(1)$ and $SU(2)$, but it is fairly clear that
it is fundamentally a discrete phenomenon rather than a continuous one. In practice, there appears to be a small experimental variation
of the Weinberg angle from a theoretical maximum of $30^\circ$.

 If we assume for the moment that the underlying discrete property is described
by an angle of exactly $30^\circ$, then this angle appears naturally as the angle between the complex numbers $i$ and $\omega$.
Since the group algebra model contains both $i$ and $\omega$, where the standard model only contains $i$, there is every prospect
that the group algebra model can explain electro-weak mixing, at least to first order. The small deviation of the Weinberg angle
from $30^\circ$ will of course require more detailed investigation.

More specifically, the finite versions of $U(1)$ and $SU(2)$ in the standard model are the scalar group of order $4$ generated by $i$, and the
quaternion group $Q_8$ generated by $i\sigma_1$ and $i\sigma_2$. But the scalars that appear in the finite model at this point are the scalars
$\omega$ and $\bar\omega$ that are used to convert from the representation $4_H$ to $4_C$, and from one complex version of $4_C$ to the other, by multiplying with the elements of order $3$ in the group $G$. This construction arises from the natural 
identification of $U(1)$ in the component $\CC$ of the group algebra with $U(1)_Y$, so that the latter is based on the finite group $Z_3$
rather than $Z_4$.

Thus the pair of Weyl spinors defined by $i$ and $-i$ in the standard model becomes a
triplet of Weyl spinors defined by $1$, $\omega$ and $\bar\omega$.
Indeed, the tensor product 
\begin{align}
2 \otimes_\CC 4_H = 4_C
\end{align}
converts the real (Majorana) spinor $4_H$ into a complex (Dirac) spinor,
and identifies the $Z_3$ copy of $U(1)$ with a $Z_4$ copy inside $Q_8$ inside $SU(2)$.

We now have three distinct angles in the group algebra model that appear to correspond in different ways to the
electro-weak mixing angle in the Standard Model. Two of them are exact values of $30^\circ$ and approximately
$28.15497^\circ$, while the third appears to describe a very small general relativistic correction to the second.
Now the value of $28.15497^\circ$ is obtained by considering only the quantum numbers of charge, weak hypercharge and
weak isospin, so has no contribution from energy. On the other hand, the value of $30^\circ$ is obtained by mapping
onto the Dirac spinor, and therefore involves both energy and mass. 

It follows that an experimental `running' of the mixing angle with the energy scale can be obtained by 
postulating that the single angle measured experimentally is in fact a weighted average of these two distinct theoretical angles.
We should therefore expect experimental values to lie in the range between $28.15497^\circ$, the theoretical limit at zero energy,
and $30^\circ$, the theoretical limit at infinite energy. Experimental values in the literature vary between about $28.7^\circ$ and
$29.3^\circ\pm.1^\circ$.

\subsection{The Dirac algebra} 
I have shown how the fermionic part of the algebra, $\HH+M_2(\CC)$, contains all the structure of the Dirac spinor, the Dirac equation,
the Dirac algebra and the basic principles of electro-weak mixing. These parts of the standard model use essentially all of the structure of this
part of the group algebra. The rest of the standard model must therefore lie in the bosonic part of the algebra $\RR + \CC + M_3(\RR)$.
I have suggested ways of constructing this part of the algebra via a finite analogue of the construction of the Dirac algebra.

As far as the Lie group $SL(2,\CC)$ is concerned, the Dirac algebra is the complex tensor product of two copies of the Dirac spinor representation.
But if this structure arises from the structure of the finite group, we must also be able to construct the Dirac algebra from the (complex) tensor product of
two copies of the $4_C$ representation of $G$. Moreover, the complex structure adds nothing of physical significance to the algebra, so that we might
as well work with the real tensor product, which has structure 1+1+2+3+3+3+3. This contains the whole of the bosonic part of the algebra, together
with a spare copy of $1+3$.

There are two ways of looking at this tensor product. One way is to look at the finite group acting on the spinors on both sides. This gives a
description of the internal symmetries of the elementary particles, without any Lie groups and therefore without any measurements or observations.
The other way is to look at the finite group acting on the vectors on the right. This gives a description of things that can be measured
by experiments that act on the vectors on the left. 

In particular, the spare copy of $1+3$ must represent things that cannot be measured, for example colours,
or direction of spin. The corresponding part of the Dirac algebra
consists of two $4$-vectors, one of which is usually interpreted as $4$-momentum, so that we might interpret the other as spacetime position,
and allocate the non-measurability to the Heisenberg uncertainty principle. That leaves us with $1+2+3+3+3$ things that can be
measured, with a macroscopic group $U(1)$ acting on the $2$, and $SL_3(\RR)$ acting on $3+3+3$. One way of obtaining the 
standard model $SU(3)$ from this structure is for both $U(1)$ and $SL(3,\RR)$ to act on $2\otimes 3=3+3$, so that we can take
complex linear combinations of elements of the Lie algebra of $SL(3,\RR)$ to make elements of the Lie algebra of
$SU(3)$. This mathematical construction is closely analogous to  that used in electro-weak theory to convert from $U(1)\times SU(2)$ to $SL(2,\RR)$
in order to built the ladder operators.

Now in removing the fourth copy of $3$ we may have removed either momentum or position, but not both. Since position is irrelevant in the
standard model, we have presumably kept the momentum. Moreover, we have a real scalar in the remaining copy of $1$, which must similarly
be the energy. We then have a macroscopic group $GL(3,\RR)$ acting on the momentum. This is a little disconcerting, since our everyday
experience is that momentum has only an $SO(3)$ symmetry. But we know from the theory of special relativity \cite{SR} that momentum mixes
with mass, and we know from the weak interaction that mass can be converted into momentum, so we should not be too surprised. 

In any case, $3+3+3$ contains a discrete set of $9$ `things' to measure, and a continuous group of dimension $9$ to do the measuring with.
Hence we can measure $9$ masses, for example of $3$ generations of electrons and up and down quarks. What we observe is that
three of these masses, the electron masses, are masses of particles whose momentum we can also measure, and they are therefore well-defined
and precisely measured masses. The other six, the quark masses, are masses of particles whose momentum we cannot measure, so that the
masses are quite ill-defined, and vary from one experiment to another. Most of all, they vary from one type of experiment to another, for
example between baryon experiments and meson experiments.

These $8$ mass ratios are closely related in the group algebra to $8$ parameters in the subalgebra $M_2(\CC)$, but the standard model
does not include a finite group relationship between these two parts of the algebra, so has an entirely different set of $8$ parameters here.
These are conventionally written in the symmetric square representation as $3\times 3$ complex matrices, and are the $8$ mixing angles in
the Cabibbo--Kobayashi--Maskawa (CKM) matrix \cite{Cabibbo,KM} and the Pontercorvo--Maki--Nakagawa--Sakata (PMNS) matrix
\cite{Pontecorvo,MNS}.

In this section I have given only a very brief sketch of how the differing masses of the three generations of fermions can be distinguished,
but without offering any physical principles by which the mass is generated. For this, I rely for now on the standard model, and defer the
question of how the experiments give the particular mass values that they do. Answering this question
requires a much deeper analysis of
the quantum structure of the experimental apparatus itself, and 
is beyond the scope of this paper. In subsequent papers \cite{octahedral,gl23model} I have attempted to
analyse the properties of mass in more detail. 

\section{Field theories}
\label{physics}
\subsection{Principles}
\label{principles}
There are at least three 
different ways of looking at the groups acting on the group algebra. One can look at the finite groups acting on both left and right, or the Lie groups acting on both left and right, or one of each (here, I assume the Lie group acts on the left, the finite group on the right). 
The general principle must surely 
be that if the Lie groups act on both sides, then we see no quantisation,
and we must recover a reasonable approximation to a description of classical physics. 
At the other extreme, if the finite group acts on both sides, then all we see are
quantum numbers, and no macroscopic variables like position, momentum, energy, mass and so on. 
Hence there are only internal variables, and nothing that can actually be measured (or `observed') from outside.

In between, we see a right-multiplication by a finite group, which flips the quantum numbers, and a left-multiplication by some Lie groups, that permits us to measure the properties of mass, momentum, energy, and so on, that are associated to particular sets of quantum numbers in particular types of experiment 
that we might contemplate performing.
This, at least, must be the general set-up in any hypothetical discrete theory of the type envisaged by Einstein \cite{Einstein1935,Einstein1954}. 

The standard model
interprets things differently, and regards the Lie groups as being intrinsic to the elementary particles, rather than
contextual properties of their embedding in spacetime.
It is clear, of course, that the calculations that are done in the standard model are essentially correct. It is only the interpretation
that must be somewhat different in a discrete model.
The discrete model must 
in fact be capable of explaining not only the standard model, but also classical physics. This is a major challenge for
any model, and success is hardly to be expected. Nevertheless, let us see how the proposed finite model deals with this challenge.
The first and most basic issue is to determine how macroscopic spacetime emerges from the discrete properties of elementary particles
and their interactions. This must be some kind of generalisation of the toy example already given in Section~\ref{analysis}.

\subsection{Classical forces and relativity}
\label{classical}
Looking at the real group algebra first from a macroscopic perspective, with the 
matrix groups acting on themselves by conjugation,
we see three real and two complex scalars, that act trivially, plus three non-trivial symmetry groups: 
\begin{eqnarray} 
SO(3) \times SO(3,1) \times SL(3,\RR).
\end{eqnarray}
Therefore, in addition to the magneto-weak `spacetime' $\HH$ with symmetry group $SO(3)$, discussed in Section~\ref{analysis},
we see an electromagnetic (special
relativistic) `spacetime' with symmetry group $SO(3,1)$, and a further space (with or without time) with symmetry group
$SL(3,\RR)$. 
The issue then is to decide how the macroscopic spacetime that underlies classical physics relates to these three different
versions of 
spacetime 
that seem to emerge in some way from different parts of particle physics. 

There is nothing in classical physics directly corresponding to the group $SL(3,\RR)$, which allows
stretching and shearing of space. But 
combining it with the Lorentz group $SO(3,1)$ acting on the same spacetime gives the group $SL(4,\RR)$
of all unimodular coordinate changes, that describes the general covariance of general relativity \cite{GR1,GR2,GR3}
as a
theory of gravity. So $SL(3,\RR)$ must in some sense describe a `fluid' gravitational `space'. 

Classically, of course, there is only one spacetime, with symmetry group $SO(3)$ 
defined by the observer. The groups $SO(3,1)$ and $SL(3,\RR)$ are then interpreted, not as different types of spacetime, but as forces. 
Clearly $SO(3,1)$ is the symmetry group of electromagnetism, as most elegantly expressed by Einstein's interpretation of Maxwell's equations
in the theory of Special Relativity.
To be more precise, the electromagnetic field is an association of a (trace $0$) element of $M_2(\CC)$ with each point in spacetime $\HH$.
It is therefore expressed mathematically as a function 
\begin{eqnarray}
\label{field}
f&:& \HH \rightarrow M_2(\CC).
\end{eqnarray}

Maxwell's equations 
then describe how the values of the electromagnetic field, in the adjoint representation of $SO(3,1)$, 
together with the inertial mass and charge, 
in the scalar part of $M_2(\CC)$, relate to the ambient physical space, in the adjoint representation of $SO(3)$,
and to time, in the scalar part of $\HH$. 
The theory of special relativity explains how the equations are invariant under coordinate changes on spacetime described by the Lorentz group $SO(3,1)$. That is, any change in spacetime coordinates defined by an element of $SO(3,1)$ can be compensated for by the
corresponding change in coordinates for the values of the electromagnetic field, by the corresponding element of $SL(2,\CC)$ acting
by conjugation on the trace $0$ matrices.

The Lorentz group can in turn be interpreted as a change of spacetime coordinates between different observers. But the model under discussion is a model appropriate to a single observer, with a fixed notion of space and time. Therefore the appropriate interpretation of the Lorentz group in this context is as a gauge group, that can be used to choose coordinates for the (electromagnetic) gauge field. The content of the theory of special relativity is then that the theory of electromagnetism does not depend on the choice of coordinates. 

\subsection{Properties of gravity}
Analogously, we must seek a gravitational field in the adjoint representation of $SL(3,\RR)$. From the point of view of an individual observer,
with $SO(3)$ symmetry, this representation splits up as the sum of a spin $1$ field, consisting of the
anti-symmetric matrices,  and a spin $2$ field, consisting of the symmetric matrices. Newton's universal theory of gravitation 
already includes both: the spin $1$ field is the gravitational field, and the spin $2$ field describes the tidal forces that arise from rotations of matter within a gravitational field (or, equivalently, rotations of the gravitational field around the matter). 

Newton's theory is very good indeed in most circumstances, but has two main drawbacks that have become apparent in the ensuing centuries. One is that it does not distinguish between the gravitational force due to a rotating body, and that due to a stationary body. The other is that it has a static gravitational field, rather than a dynamic field that propagates at a finite speed (presumably the speed of light).

The former drawback is (at least partly) addressed by Einstein's theory of general relativity, in the sense that
rotations of the observer (or test particle) are taken into account. But it is not clear that rotations of the gravitating body are fully
accounted for, 
nor 
the finite speed of propagation.
It is therefore not required for a new model to reproduce general relativity exactly, but merely to reproduce general relativity in the limit where the effects of both rotation and the finite speed of light can be ignored. 

There are indeed 
two circumstances in which such effects might already have been observed in practice. One is in the effects on fast-moving satellites as they rotate very fast around the Earth, and
% The former 
is known as the flyby anomaly \cite{flyby}. It is claimed by Hafele \cite{Hafele} that this anomaly can be entirely explained
 by the finite speed of propagation of the gravitational field of a massive rotating object. The effect is observed as a consequence of the
 interplay between the rotation and the speed of propagation, rather than one or other in isolation.
 
 The other is at the edges of galaxies, in which the gravitational attraction of the galactic centre can take hundreds of thousands of years to reach the outer edges. 
 Many proposals
 have been made for theoretical explanations for the observed effects, most notably the hypotheses of (a) dark matter, or (b) modified Newtonian 
 dynamics, or MOND \cite{Milgrom1,Milgrom2,Milgrom3,MOND,MOND2}. 
 Indeed, Yahalom \cite{Yahalom} proposes an explanation of these astronomical observations based entirely on %that takes into account either 
 the finite speed of propagation
 of the gravitational field.
 In this work, it appears to be unnecessary to take into account the very fast rotation of the super-massive black hole (or other objects) at the centre of the galaxy. 

The proposed finite model permits a separation of the $SL(3,\RR)$ symmetries of gravitational space from the $SO(3,1)$ symmetries of
electromagnetic spacetime, and therefore permits a macroscopic theory of gravity that is independent of time, and therefore independent
of the finite speed of light. But to do that it requires the symmetry group of the gravitating body to be extended from $SO(3)$ to $SL(3,\RR)$.
In other words, the centre of the galaxy must be regarded as a rapidly rotating fluid rather than  a solid object. Given that a typical galaxy has many billions of stars,
this is surely a reasonable hypothesis to make. 

The general principle of relativity says that it makes no difference
to the theory whether the group $SL(3,\RR)$ is attached to the galactic centre or to the periphery, so that we can interpret the theory
either way. 
The latter interpretation seems to be closely related to the dark matter hypothesis: that is the galactic centre is regarded as a Newtonian
point mass, and all the otherwise unexplained behaviour of the outer stars must then be due to some `dark matter halo' surrounding the galaxy. 
The former interpretation is more closely related to the MOND hypothesis: that there is some new unexplained gravitational force,
that 
can only be observed when the Newtonian gravitational field is extremely weak. 

The proposed finite model
does not favour the dark matter interpretation, 
since it contains no 
type of matter that is not already in the SMPP. 
In the MOND context, it suggests (very speculatively) that the new force might be explained in terms of the very fast rotation of the extremely massive
galactic centre.
Most importantly, the group $SL(3,\RR)$ is now interpreted as the gauge group of the theory of gravity, and is completely separate
from the gauge group $SO(3,1)$ of classical electromagnetism.
This might explain why attempts to quantise general relativity with
a gauge group $SL(4,\RR)$ or $GL(4,\RR)$ are regarded as having failed \cite{GL4R1,GL4R2}.

On the other hand, there is no room in the group algebra model for new forces, so the only realistic possibility is 
a re-interpretation of old forces. In particular, the group $GL(4,\RR)$ of general covariance does not appear explicitly in the model,
which contains only $GL(3,\RR)$ acting on coordinates of space. There is then a scalar (invariant) gravitational time, in contrast to the
covariant electromagnetic spacetime. The best way to reconcile these two inconsistent definitions of time is to put the
time delay of propagating forces explicitly into the model of gravity.

\subsection{Emergence of classical physics}
The great virtue of working with a group algebra is that one can switch viewpoint from the algebra to the group and back again, and therefore
combine the continuous symmetries of classical physics and relativity with the discrete symmetries of quantum mechanics.
The algebra can be acted on by both left multiplication and by right multiplication, either by the whole algebra or just by the finite group.
The linearity of quantum mechanics is obtained by letting the continuous groups act on one side and the finite group on the other.

This allows us to split up the classical fields into the rows of the matrices, which gives us three classical vector fields, for Newtonian gravity
and electromagnetism, together with a charge in $\RR$ and a mass plane in $\CC$. The gravitational part of this therefore consists of the five
real dimensions of
$\CC+\RR^3$, while the elecromagnetic part consists of the seven real dimensions of $\RR+\CC^3$.
In General Relativity, however, the continuous symmetries of spacetime act on both rows and columns, so that the gravitational part of
the theory lies in the symmetric matrices, and the electromagnetic part lies in the anti-symmetric matrices,
which split under the action of the finite group as
\begin{align}
S^2(\RR+\RR^3) & = 2\RR+\CC+2\RR^3\cr
\Lambda^2(\RR+\RR^3) & = \CC^3.
\end{align}
The occurrence of two copies of $\RR+\RR^3$ in the symmetric matrices might indicate that there is some redundancy in General Relativity,
and that for the purposes of quantisation we only need to consider the spatial part of the symmetric tensor
\begin{align}
S^2(\RR^3) = \RR + \CC + \RR^3.
\end{align}

Under the action of $SL(3,\RR)$, this representation is irreducible, while as a representation of $SO(3)$ it splits as $1+5$, where $5$
represents a spin $2$ representation on the space $\CC+\RR^3$.
This is the fundamental reason why quantum gravity is expected to contain a spin $2$ graviton.
But in the model as described, there is no such graviton, simply a splitting into the classical Newtonian field, plus a mass plane.
General Relativity now has the effect of moving the scalars $\RR$ from electromagnetism, where it represents charge,
to gravity, where it represents the scalar curvature of space. 

Indeed, the whole of the 12-dimensional group algebra splits up into 6 dimensions of gravity plus 6 dimensions of 
electromagnetism. 
Unification of gravity with electromagnetism can then be expressed in the tensor product
\begin{align}
(\RR+\RR^3) \otimes \RR^3 = \RR+\CC+\RR^3+\CC^3,
\end{align}
in which $\RR+\CC+\RR^3$ represents Einsteinian gravity, and $\CC^3$ represents Einsteinian electromagnetism.
This unified force can then be quantised by expressing the first factor as the (quaternionic) square of weak hypercharge and weak isospin, to give
a unified quantum field in the form
\begin{align}
\HH \otimes _\HH \HH \otimes _\RR \RR^3.
\end{align}

\subsection{Emergence of quantum gravity}
In particular, the model quantises gravity using the representation $\RR+\CC+\RR^3$ of the finite group. In effect, $\RR^3$
represents the Newtonian gravitational field, which propagates at the speed of light in the same way that the electromagnetic
field $\CC^3$ does. The other half $\RR+\CC$ represents the weak force, with its broken symmetry, to create the mass and charge.
In other words, General Relativity mixes Newtonian gravity with the weak force. 

This explains why the electroweak mixing angle comes in two different forms, one that is fixed, and does not involve gravity,
and one that depends on the masses of elementary particles measured in a particular gravitational field. In effect there is a `mixing angle'
between gravity and the weak force, that is the difference between these two versions of the electro-weak mixing angle.
Since this angle is very small, slightly less than $38''$, it is undetectable within the noise of particle physics experiments.
But now that we have a mathematical model that predicts this discrepancy, we can look for a physical explanation.

Any quantisation of gravity must ultimately rely on a massless particle (i.e. graviton) that propagates the force across space.
The speed of such propagation must under any reasonable assumptions about the nature of spacetime be equal to the speed of light.
 The mixing of $\CC$ with $\RR^3$ is, almost by definition, a tidal effect, in which the main players are the Sun and the Moon.
 The model does not allow us to use the time coordinate, and therefore does not allow us to use the magnitude of the gravitational forces.
 What it does allow us to use are angles in space.
 
 The crucial property of a long-range force that propagates at a finite speed is that Newton's third law does not hold \cite{Hafele,Yahalom}: 
 in the time it takes for the two bodies to exchange a graviton, they have moved significantly from their previous positions.
 The critical parameter, therefore, is the angle between the action and the reaction. It is easy to calculate this angle in the case of the
 tides caused by the Sun and the Moon.  The Moon is close by, and the angle by which the Moon has moved in its orbit
 during the time that a signal at the speed of light travels from the Earth to the Moon and back (around 2.7 seconds), is 
 not much more than $1''$. The corresponding angle for the Sun is around $41''$.
 
 This effect would appear, therefore, to be the main contributor to the discrepancy between the two calculated values of the
 weak mixing angle. There are obviously many other contaminating effects, not all of which we can examine here. 
 But if a change in the direction of gravity between the active and passive phases is important, then it is reasonable to suppose that
 a change in the direction of gravity between different parts of the experiment is also important. Many particle physics 
 experiments are very large, of the order of 5km in diameter, so that the direction of the gravitational field changes by something 
 of the order of $3'$ in such an experiment.
 It is therefore clear that contamination from quantum gravity can cause systematic errors of this magnitude in any 
 measurement of parameters of the weak force. The CODATA quoted standard uncertainty of the Weinberg angle of about $1.3'$
 is therefore almost certainly an underestimate. 

\section{Waves and particles} 
\label{waveparticle}
\subsection{Quantum field theory}
The above discussion suggests interpreting the $24$ dimensions of the group algebra as classical fields, divided into $7$ dimensions of scalar fields,
$3$ dimensions of space, $6$ dimensions of electromagnetic field, and $8$ dimensions of gravitational field. These fields are mathematical
constructs, rather than physically real objects. For example, we do not want to interpret the $3$ dimensions of the `space' field as an `aether'
in the 19th century sense. What is physically real, however, is the propagation of the fields through space, or to be more precise, the propagation
in spacetime. There is also a real physical distinction between the $12$ fermionic `matter' fields and the $12$ bosonic `force' fields.

I have not precisely identified all of the scalar fields, though four of them appear to be time, inertial mass, charge and gravitational mass.
Energy would seem to be separate from both types of mass, and there appears to be a second type of (neutral) charge. This leaves one more
scalar, which might be identified with the Higgs field in the standard model. Whatever identification is eventually decided on, these scalar
fields propagate by being attached to matter particles. In other words, it is the movement of matter that defines the propagation of the
scalar fields. 

Just as in the classical case, the fields must be modelled as functions from spacetime to the appropriate subspace of the group algebra.
In the quantum case, this function decomposes as a sum of `infinitesimal' quanta, that individually are linear functions. These linear functions
themselves lie in representations of the finite group $G$, so that we can use the representation theory to analyse the structure of the
quantum fields.
Since the spacetime representation $\HH$ is self-dual, the quantum fields can be regarded as lying in the tensor product of
$\HH$ with the appropriate representation. 

These tensor products then describe the particles that mediate the corresponding forces.   Since the spacetime 
 representation is fermionic in this model, the bosonic fields are carried by fermionic
particles, and the fermionic fields are carried by bosonic particles. In particular, four of the scalars are bosonic, so carried by fermions, and
three of the scalars are fermionic, so carried by bosons.

The electromagnetic field values are fermionic in this model, and the field 
is therefore propagated by bosons. This agrees with the standard model, in which the propagation is effected by photons, which are parametrised by $3$ dimensions of momentum in $2$ distinct polarisations.
The three fermionic scalars, which might be interpreted as inertial mass (dual to Euclidean time, as opposed to energy,
which is dual to Lorentzian time), electric and neutral charge, are 
carried by the weak bosons, that is the $Z$ and $W$ bosons.
Conversely, in this model the gravitational field values are bosonic, and the field is therefore propagated by fermions. 

No such process exists in the standard model,
in which all fields are assumed to be carried by bosons.
But it would appear to be consistent with experiment 
to suppose that these propagators are neutrinos,
parametrised by $3$ dimensions of momentum in $3$ distinct generations. This parametrisation would result in a $9$-dimensional
gravitational field, however, and the model (as well as observation) supports only an $8$-dimensional field. 
In other words, the $3$ generations of neutrinos cannot be linearly independent. This means that macroscopic rotations act not only on
the momentum coordinates, but also on the generation coordinates. The model therefore predicts that the generation of a neutrino
is not an invariant. Indeed,
the non-invariance of neutrino generation is well-attested experimentally, and goes by the name of neutrino oscillation
\cite{oscillation,neutrinos,SNO}.

At this point, the neutrinos appear to have taken over the group $SL(3,\RR)$, that was originally supposed to be allocated to the strong force,
and in particular to the $8$ gluons. Indeed, since the group $SL(3,\RR)$ in the proposed model acts on a bosonic field, the corresponding
mediators must be fermions. We can reconcile the two viewpoints by interpreting the gluons as representing the \emph{values} of the quantum field,
rather than the field itself, which is a \emph{function}. Then the mediators can be interpreted as virtual neutrinos, and the gluons as pairs of virtual neutrinos. However, the model suggests that it may be better not to interpret the gluons as particles at all, but only as symmetries.

At the same time, we need to address the distinction between the group $SL(3,\RR)$ used here and the group $SU(3)$ used in the
standard model. The latter is a compact group, and fixes a complex inner product, so describes rigid symmetries of a complex $3$-space.
The phenomenon of asymptotic freedom \cite{freedom1,freedom2} suggests that such rigidity does not in fact characterise the strong force. The use of
$SL(3,\RR)$, on the other hand, suggests complete freedom to change scale in one direction, providing this is compensated for in
another direction. In other words, replacing the (confined) gluons by the (free) neutrinos seems to require a split group, just as the (free)
photons are described by the split group
$SL(2,\CC)$.

\subsection
{Elementary particles and the standard model}
Let us now turn our attention to the mixed case, with the finite group acting on the right and the Lie groups acting on the left. This is the domain
of the standard model, where measurements of macroscopic variables such as mass, momentum, energy, angular momentum,
magnetic moments and so on are made on individual quanta or `elementary particles'. 

We then have 24 discrete objects on the right
that we can measure, and 24 degrees of freedom for the operators on the left that define what we can measure. 
If we make a specific choice of the $24$ independent things we want to measure, we obtain $24$ (dimensionless) real numbers that describe everything
there is to know about how the elementary particles behave. 
The standard model has made a particular choice, and has measured precisely $24$ independent things. These $24$ things are usually described as
$12$ fermion masses ($3$ generations each of neutrino, electron, up and down quark), $3$ boson masses (the $Z$, $W$ and Higgs bosons), (hence $14$ mass ratios), $2$ coupling constants (the fine-structure constant and the strong coupling constant) and $8$ mixing angles ($4$ each in the
CKM and PMNS matrices)
\cite{Cabibbo,KM,Pontecorvo,MNS}, and are regarded as the fundamental parameters.

The standard model therefore has exactly the right number of parameters to contain a complete description of quantum reality. It produces the right answers, because it contains enough variables and enough equations to calculate everything that can be calculated. There is nothing wrong
with the standard model. It is a mathematically correct, and complete, model of everything that can possibly happen. So 
why are people
still looking for a `theory of everything', if the standard model already \emph{is} a theory of everything?

The main reason is that we don't understand where the $24$ parameters come from. We understand perfectly well that the Lorentz group
can be interpreted either as a change of coordinates on spacetime between two observers, or as a gauge group for electromagnetism, and that
these two interpretations are equivalent, 
so that the theory of electromagnetism is the same for all observers. We understand perhaps a little less well that the same applies to $SL(3,\RR)$ and the theory of gravity. 

So why don't we understand that the same applies to the
gauge groups of the standard model? The groups describe the relationship between the observer and the observed, and we have a choice
between regarding the groups as acting on the observed (as in the standard model) or the observer (as in relativity). 
The incompatibility of the standard model with general relativity may be nothing more than the fact that the two
disciplines have made incompatible choices of interpretation.
So to resolve the issue we 
need to make a consistent choice of interpretation. For practical purposes, it makes no difference which choice we make.
Either way we do the same calculations, get the same answers, and reconcile them with experiment. 

But philosophically, there is no contest. 
The general principle of relativity is such a powerful and obvious philosophical principle that we should not under any circumstances contemplate abandoning it. 
The consequence of this philosophical viewpoint is that the gauge groups of the standard model have to be transferred to act on the observer, so that
(most of) the $24$ parameters become parameters associated with the experiment, the environment and the observer, rather than with the
elementary particles themselves. 

I re-iterate that both points of view are mathematically valid, but only one of them is philosophically valid.
Moreover, since many of the $24$ parameters are known to vary with the experiment, and in particular to `run' with the energy scale, 
it is also the case that only one of the two viewpoints is
physically valid. It is not physically reasonable to treat a parameter as a universal constant,
if experiment shows that it is not.
In the model I am proposing, $7$ of the $24$ parameters are scalars, and can therefore be taken as universal constants.
The other $17$ must be regarded as properties of the experiment. 

If we look carefully at the experimental evidence, $6$ or $7$ of the $14$ mass
ratios  
show evidence of not being constant, namely the $6$ quark masses and the $Z/W$ mass ratio. All $8$ of the mixing angles are similarly suspect.
Finally, there is no positive, model-independent, evidence of different masses for the three generations of neutrinos, or indeed any mass different from zero. Moreover, the phenomenon of
neutrino oscillation suggests that there is no intrinsic difference between the three generations of neutrino anyway, which add another $2$ or $3$
parameters to make up a total of between $16$ and $18$.  

\subsection{The wave-function}
The one thing that one might really hope for from a finite model of elementary particles is an insight into the measurement problem.
This was really the focus of Einstein's objections to quantum mechanics throughout his life, and although his attacks, most notably
the EPR paradox \cite{EPR}, were never enough to sink the ship, the fundamental problem has not gone away.
To simplify the problem almost to the point of caricature, we may ask, what \emph{is} the wave-function, and how does it `collapse'?

In the group algebra model, the wave-function is implemented at the middle level, with the finite group acting on the right
and the Lie groups on the left. 
The typical example is equation (\ref{field}), which describes a
function from spacetime to the Dirac spinors. The `collapse' is some kind of operation that moves down to the discrete level,
with the finite group on both sides. In practical terms, the wave-function describes the quantum field in the experiment, and the collapse describes
the result of the experiment.

From this point of view, the measurement problem arises from the assumption that the wave-function is intrinsic to the elementary
particle under investigation. The finite group model, however, does not allow any continuous variables to be associated with an elementary
particle. The continuous variables are always associated with the macroscopic measuring apparatus. The model, in other words, does not
\emph{solve} the measurement problem. It \emph{does not find} the measurement problem.
Philosophically, this is always the best way to solve problems, namely, to realise that, if looked at in the right way, the problem does not exist.
It should be noted, however, that this proposed change in interpretation has no effect at all on the calculations of probabilities and other
continuous variables in quantum mechanics and the standard model.

\subsection{The real universe}
The structure of the model indicates that of the $24$ unexplained dimensionless parameters of the standard model, exactly $7$
are universal constants, and the rest are dependent in some way on the experiment, the environment and/or the observer. There is nothing
like enough detail in the model as so far developed to indicate which parameters should be regarded as universal constants, nor how the
other parameters vary. Indeed, there may be 
a certain amount of choice as to which parameters can be defined as constant, so that the other
parameters can be calibrated against them. 

What is clear, however, is that we must be prepared for the possibility that
certain parameters that we are certain
must be universal constants, may not be so.
At this stage, we can do little more than speculate on these matters, and apply some educated guesswork. There is no real need to use
the same fundamental parameters as the standard model, but other mass ratios such as those between electron, proton and neutron might
also be worth looking at. The model suggests that approximately $8$ of the fundamental parameters are dependent in some way on the
gravitational field, including tidal forces. 

At first glance, the very idea is preposterous. But on closer inspection, one realises that while the
standard model is in theory based in an inertial frame, the experiments that measure the fundamental parameters are not done in an
inertial frame, but in a frame which moves with the Earth. It is therefore very hard to rule out, on experimental grounds, the possibility
that some of the measurements depend crucially on some dimensionless property of the tides that is constant over all experiments done
on the Earth.

There are four basic dimensionless parameters of the tides on the Earth, of which two are obvious spatial angles: the angle of tilt of the Earth's axis,
and the inclination of the Moon's orbit to the ecliptic. The other two are the ratios of the periods of rotation/revolution, so are angles in Euclidean spacetime. Of course, these
parameters are not precisely constant, and there are other effects that might be expected to contaminate the results, such as the
gravitational pull of Jupiter. It will therefore be somewhat tricky to distinguish a real correlation from an unwanted coincidence.

Some $8$ such correlations/coincidences were presented in \cite{INI13}, without any very solid justification, but with the suggestion
that, while some of them may indeed be pure coincidences, it is very unlikely that they are all coincidences. The parameters discussed
in that paper include the electron/proton/neutron mass ratios, the pion mass ratio, the kaon mass ratio, the kaon/eta mass ratio,
the Cabibbo angle, the Weinberg angle and the CP-violating phase in the CKM matrix. In addition, the paper discusses a number of other
equations that seem to suggest that some of the fundamental parameters may not be independent of each other. 

This may be a
somewhat less unpalatable suggestion than the suggestion, made originally by Einstein \cite{gravimass} more than a century ago, 
that they might depend on the gravitational field! On the other hand, the paper \cite{INI13} also shows that the observed CP-violating
behaviour of neutral kaons \cite{CP} is quantitatively consistent with the hypothesis that the effect is caused by the small difference in the
direction of the gravitational field between the two ends of the experiment.  Further analysis of these various questions in presented
in \cite{gl23model}, where more detailed proposals are made for relating specific parameters to specific elements of the model.

\subsection{Composite particles}
One of the early successes of Gell-Mann's eightfold way was the explanation of the meson octet/nonet and the baryon octet and
decuplet in terms of representations of $SU(3)$, and in particular the prediction of the $\Omega^-$ particle to complete the baryon
decuplet. These are all particles that in the standard model are made up of up, down and strange quarks, and their anti-quarks.
The standard way to break the symmetry of these quarks is by using strangeness as a proxy for mass, so that the adjoint representation
of $SU(3)$ splits as 3+(2+2)+1, interpreted as three pions, four kaons and the eta meson.

But strangeness is not a conserved quantity, and the finite model suggests a splitting instead in terms of charge, which is conserved.
Then the two pairs can be interpreted as charged pairs, of pions and kaons, and the singlet as a neutral pion. These three sets are, experimentally,
mass eigenstates.
 This leaves the triplet
for two neutral kaons and the eta meson, with two different masses. Or, as suggested in \cite{INI13}, we could keep the mass eigenstate property,
and allocate three neutral kaons, so that it becomes unnecessary to use quantum superpositions to explain the properties
of neutral kaons. Indeed, the finite model requires this, since it permits quantum superposition only as a property of the experiment,
not as an intrinsic property of a physical particle.

I have already mentioned that the finite model describes the neutral pion as $u\bar u$ rather than a quantum superposition of $u\bar u$
and $d\bar d$. Similarly, the neutral kaons must then be allocated to $d\bar d$, $d\bar s$ and $s\bar d$, and the eta meson to $s\bar s$.
This allocation is necessarily slightly different from the standard model, in order to avoid quantum superpositions.
All of this analysis is based on the tensor product representation
\begin{eqnarray}
(1+2)\otimes (1+2) &=& 1+1+2+2+3
\end{eqnarray}
which applies equally to the finite group $G$ and the Lie group $SU(2)$.

Another possible interpretation of this $9$-dimensional representation is as the $9$ fundamental particles of matter: the
three generations of electrons, up and down quarks. In this interpretation the irreducibles represent charge eigenstates,
with a further splitting of the quarks into $3$ light and $3$ heavy quarks. Thus the two singlets are the up and bottom quarks,
and the two pairs are down/strange and charm/top. In particular, in this model the three generations of quarks look quite different from the
three generations of leptons.

Baryons are then described in the standard model by parts of the tensor product of three copies of $1+2$
(or $1+1+2+2$). The full structure is
\begin{eqnarray} (1+2)^3 &=&
(1+2)\otimes (1+1+2+2+3)\cr
&=&
1+1+1+1+1+2+2+2+2+2+3+3+3+3\cr
&=& 1 + (1+2+2+3) + (1+2+2+3) + (1+1+2+3+3)
\end{eqnarray}
where the parentheses indicate the splitting for $SU(3)$ as $1+8+8+10$. The splitting of the baryon octet on the basis of mass as one $\Lambda$, two nucleons,
two $\Xi$ and three $\Sigma$ baryons is clear enough, but they can also be split on the basis of charge if so desired. The baryon decuplet splits as
$1+2+3+4$ either way, but there is no clear reason why the $4$ should split as $1+3$.

The finite model contains a number of other related representations that might provide extra insight. For example, replacing the broken $1+2$ symmetry
by the unbroken $3$ symmetry gives us the representation
\begin{eqnarray}
(1+2)\otimes 3 \cong 3\otimes 3 &\cong& 1+2+3+3
\end{eqnarray}
which offers a splitting of the baryon octet into a pair of $uds$ baryons (the $\Lambda$ and $\Sigma^0$) and two triplets,
one of $p$, $\Sigma^-$ and $\Xi^0$, the other of $n$, $\Sigma^+$ and $\Xi^-$. Each of these triplets contains one particle with each 
charge $+$, $-$ and $0$, and as Coleman and Glashow observed in the early 1960s, the two triplets have the same total mass
\cite{ColemanGlashow}. Alternatively, two copies of the same representation $2+3+3$ arise from $4_H\otimes 4_C$.

It is noticeable that the exact Coleman--Glashow relation applies to triplets of particles, while no such exact relations appears to exist
for doublets (such as proton/neutron). An exact formula for the generation triplet of electron, muon and tau particle, each paired with
a proton to cancel the charge, was presented in \cite{INI13}: but \emph{why} should the total mass of these six particles be equal to the
mass of five neutrons? A possible reason for this mass formula is proposed in \cite{gl23model}, where the concept of spin is
analysed in detail, with a total of $15$ different types of spin, that can be combined into particles in various different ways.
This reduction from the standard model count of $24$ (four colours, three generations, left and right) is obtained by not allowing the
quarks to have independent left and right spins.

 \subsection{Differential equations}
The equations of relativity arise from decompositions of the `bosonic' part of the algebra, in the form
\begin{align}
(\RR+\RR^3) \otimes \RR^3 & = \RR+\CC+\RR^3+\CC^3\cr
& = (\RR+\RR^3) \otimes (\RR+\CC)
\end{align}
There are of course $12$ independent equations here, which can be split into six Maxwell equations for the
electromagnetic field $\CC^3$ plus six Einstein equations for the gravitational field $\RR+\CC+\RR^3$.
But there is some mixing between the two that becomes very apparent when considered from the quantum point of view,
so that it is easier to keep all terms together and describe both forces at once.

Finally, a complete unification of all the forces of nature can be obtained from expressing the full group algebra in the form
\begin{align}
\RR^3 \otimes (\RR + \RR^3 + \HH) = (\RR+\CC)\otimes(\RR+\RR^3+\HH).
\end{align}
This equation is a shorthand for a complete set of $24$ linear differential equations that combine gravity, electromagnetism, and the
nuclear forces into a single universal quantised force law. The terms in $\RR+\RR^3$ are the long-range classical forces,
and the terms in $\HH$ are the short-range quantum forces. The splitting $\RR+\CC$ separates the weak nuclear force in $\RR$ from
the strong nuclear force in $\CC$, and separates charge in $\RR$ from mass in $\CC$. But all of these splittings are local
approximations to parts of the unified theory. No piece can be completely separated form any other. 

Perhaps the most important point of all is that the short-range forces
generate long-range forces (so that particles create waves) 
via the equation
\begin{align}
\HH\otimes_\HH \HH = \RR+\RR^3.
\end{align}
Here I am using the real part of $\HH$ to denote differentiation $\partial_{ct}$ with respect to time, and the imaginary part
to denote differentiation with respect to space $(\partial_x,\partial_y,\partial_z)$. 
The quaternionic tensor product enforces a quaternion conjugation linking the 
left-multiplication and right-multiplication, so that when we differentiate on both sides simultaneously, all the cross terms disappear,
and we get the second derivatives appearing naturally in $\RR+\RR^3$. Therefore we get the standard $3$-dimensional wave equation
for gravitational waves in a vacuum.

This allows us to attach waves to particles by factorising
\begin{align}
\RR + \RR^3 + \HH & = \HH\otimes_\HH \HH + \HH\cr
&=  \HH \otimes  (\HH + \RR)
\end{align}
which can be interpreted as attaching charge in $\RR$ to mass in $\HH$. This is, ultimately, what enables us to allocate charge and four independent
masses to the proton and three generations of electron, and to use the linearity of $\HH+\RR$ to derive the mass equation (\ref{emutau})
that inspired this research.

\section{Other possibilities}
\label{other}
\subsection{The binary octahedral group}
There are one or two places in the proposed model where it may be felt that it does not closely enough resemble the standard model.
For example, the split group $SL(3,\RR)$ stands in for the compact group $SU(3)$ as a gauge group for the strong force, while
the weak gauge group $SU(2)$ and the spin group $SL(2,\CC)$ seem to have at least partly changed places, so that the
symmetry-breaking to $SL(2,\RR)$ appears somewhat confusing. If we extend from the binary tetrahedral group $2.Alt(4)$
to the binary octahedral group $2.Sym(4)$, then there is more structure available in the group algebra that might allow us to address these issues.
See \cite{octahedral,gl23model} for more details.

The character table of the binary octahedral group is as follows:
\begin{eqnarray}
\begin{array}{cccccccc}
1&-1&i&v&-v&(i+j)/\sqrt2 & (1+i)/\sqrt2 & -(1+i)/\sqrt2\cr\hline
1&1&1&1&1&1&1&1\cr
1&1&1&1&1&-1&-1&-1\cr
2&2&2&-1&-1&0&0&0\cr
3&3&-1&0&0&1&-1&-1\cr
3&3&-1&0&0&-1&1&1\cr
2&-2&0&-1&1&0&\sqrt2 & -\sqrt2\cr
2&-2&0&-1&1&0& -\sqrt2 & \sqrt 2\cr
4&-4&0&1&-1&0&0&0\cr\hline
\end{array}
\end{eqnarray}
The structure of the real group algebra is correspondingly
\begin{eqnarray}
2\RR + M_2(\RR) + 2M_3(\RR) + 2\HH + M_4(\RR).
\end{eqnarray}
There are three chiral pairs of representations, differing only by changing the sign on the odd permutations of $Sym(4)$.
Each pair can be given a natural complex structure, so that, for example, $2\RR$ becomes $\CC$, and the odd
permutations act as complex conjugation. Similarly, the $3$-dimensional representations combine to give a $3$-dimensional complex space,
on which we can put a group $SU(3)$ or $SO(3,\CC)$ depending on the interpretation we wish to give to it.
The two copies of the quaternion algebra can similarly be combined into a $4$-dimensional complex space, that would seem to
have a natural interpretation as the space of Dirac spinors, divided into left-handed and right-handed pieces.

Indeed, the standard model effectively puts two different complex structures onto the Dirac spinor. The one just described is
implemented with $i\gamma_5$, and a separate complex structure is implemented with $i$, that converts the quaternion algebra
$\HH$ into a complex $2$-space $\CC^2$. On the other hand, the complex structure on the pair of $1$-dimensional
representations is implemented with $i$, not with $i\gamma_5$. This appears to be an inconsistency %in 
with the standard model,
and is discussed in detail in \cite{Clifford}.

This leaves $M_2(\RR)$ to contain $SL(2,\RR)$, as a broken-symmetry version of the weak gauge group, and
$M_4(\RR)$ to contain $SL(4,\RR)$ or $GL(4,\RR)$ as a potential gauge group for gravity.
This allocation is rather different from the one suggested earlier, but matches the standard model rather more closely.
In particular, the bosonic part of the algebra now supports the gauge group
\begin{eqnarray}
U(1) \times SU(3) \times SL(2,\RR)
\end{eqnarray}
with compact $SU(3)$ for the strong force, and split $SL(2,\RR)$ for the weak force with broken symmetry.

\subsection{Differences from the standard model}
The fermionic part of the algebra now exhibits some interesting differences from the standard model. The spin group $SL(2,\CC)$
has been replaced on the one hand by $2\HH$, containing $SU(2)\times SU(2)$, and on the other hand by $SL(4,\RR)$.
Thus we have all the groups we need for both the standard model of particle physics and for general relativity, but they are mixed
together in unexpected ways. This may be a good thing, because it might explain why the `obvious' methods of quantising
general relativity don't work.
In particular, we see that gravity has a `fermionic' rather than bosonic gauge group, which again suggests the need for neutrinos
as mediators for quantum gravity. 

Moreover, the group that acts on the Dirac spinors in this model is not
\begin{eqnarray}
Spin(3,1)&\cong &SL(2,\CC)
\end{eqnarray}
 as in the standard model,
but 
\begin{eqnarray}
Spin(4)&\cong& SU(2)\times SU(2).
\end{eqnarray}
 Both the groups $SL(2,\CC)$ and $SO(3,1)$ are subgroups of $SL(4,\RR)$,
where they can be interpreted as acting on spacetime, either macroscopic or `internal'. But there is no mathematical or 
physical reason for regarding these two groups as being in any sense `the same'. Any similarities between them lie
at the level of analogy (as in the original work of Dirac \cite{Dirac}, and as also explained in modern textbooks \cite{Griffiths}), 
not at the level either of
mathematical equality or of physical interpretation.

One possible advantage of using the binary octahedral group rather than the binary tetrahedral group is that the full general
linear group $GL(4,\RR)$ is available for implementing a quantisation of general relativity. But this might instead be a
disadvantage. Indeed, the fermionic part of the algebra is now essentially
\begin{eqnarray}
SU(2) \times SU(2) \times SL(4,\RR)
\end{eqnarray}
which looks very similar to the gauge group of the Pati--Salam model
\begin{eqnarray}
SU(2) \times SU(2) \times SU(4).
\end{eqnarray}
The two copies of $SU(2)$ in both cases act on the left-handed and right-handed parts of the Dirac spinor, and the remaining $15$
degrees of freedom form either the split or compact real form of the Lie group of type $A_3$. 
This suggests that the new model may suffer from the same problems as the Pati--Salam model.

\subsection{The octahedral group}
If so, then perhaps we should regard the whole of the fermionic part of the algebra as being `ficititious' in the sense that it provides a
very useful calculational tool, but should not necessarily be regarded as physically real. This is certainly a valid interpretation of the
Dirac spinor, and perhaps it is also a valid interpretation of spacetime. Then the question arises, whether it is possible to describe
physics using only the bosonic part of the algebra, without the fermionic part. It seems unlikely, but perhaps it is not impossible.

The finite group then becomes the octahedral group $Sym(4)$. This group does not contain the quaternion group $Q_8$, but it does contain
the dihedral group $D_8$. Hence we would need to abandon the idea of using $Q_8$ to describe spin, and use $D_8$ instead.
Then we see that $Sym(4)$ contains three copies of $D_8$, compared to the single copy of $Q_8$ in $2.Sym(4)$.
This gives us two potentially major advantages. 
First, the broken symmetry of $D_8$
relative to the unbroken symmetry of $Q_8$ permits us to break the symmetry between leptons and baryons, or between leptons and quarks,
or between left-handed and right-handed `spins', or even between weak isospin and these two types of spin.
Second, the three copies could be used for three generations.
Now 
restricting the bosonic characters to $D_8$ we obtain 
\begin{eqnarray}
\begin{array}{ccccc}
1 & i & j &(j+k)/\sqrt 2 & (1+i)/\sqrt 2\cr\hline
1&1&1&1&1\cr
1&1&1&-1&-1\cr
2&2&2&0&0\cr
3&-1&-1&1&-1\cr
3&-1&-1&-1&1\cr\hline
\end{array}
\end{eqnarray}
compared to the table of irreducible characters
\begin{eqnarray}
\begin{array}{ccccc}
1 & i & j &(j+k)/\sqrt 2 & (1+i)/\sqrt 2\cr\hline
1a&1&1&1&1\cr
1b&1&1&-1&-1\cr
1c&1&-1&1&-1\cr
1d&1&-1&-1&1\cr
2&-2&0&0&0\cr\hline
\end{array}
\end{eqnarray}

In particular, the $2$-dimensional representation of $Sym(4)$, which describes the weak doublets, restricts as $1a+1b$, that is the
same as the restriction of $1+1$ that describes the electromagnetic $U(1)$. Hence on restriction to a single generation these two
representations are indistinguishable from each other, and can be mixed together in arbitrary proportions. Exactly as is done in
the standard model, using the Weinberg angle. The representation $1b$ contains a sign-change which we might choose to interpret
as representing weak isospin, for example.

The chiral pair of $3$-dimensional characters restrict as $2+1c$ and $2+1d$. Here the $2$-dimensional representation is the
spin representation, analogous to the $2$-dimensional spin representation of $Q_8$. The distinction between $1c$ and $1d$ allows us to
make one further distinction, which we might choose to be that between left-handed and right-handed spin, or perhaps something
more tangible. For example, we might use the four $1$-dimensional representations for the four particles neutrino, electron, proton and
neutron, as we did originally for $Q_8$, 
so that the leptons arise from $1$-dimensional representations of $Sym(4)$, and the baryons from $3$-dimensional
representations. Of course, this discussion only involves the combinatorial structure, and does not tell us anything about the continuous structure.

 \section{Conclusion}
 I have tried to show that by putting the finite symmetries of elementary particles at the forefront of the theory, rather than the
unitary gauge groups, it is possible to gain insights into the structure of the standard model that cannot be easily obtained from the
standard approaches using Lie groups, Lie algebras and Clifford algebras. In particular, by interpreting $Z_3$ as a group of symmetries of the three generations, and $Z_4$ as a group of symmetries of the four fermion types in a single generation, we obtain simple geometric derivations for the
electron-muon neutrino mixing angle and the electro-weak mixing angle. Converting these finite groups to Lie groups gives two 
distinct gauge groups isomorphic to $U(1)$, which play the roles of $U(1)_{em}$ and $U(1)_Y$ respectively.

In the case of non-abelian gauge groups, the quaternion group $Q_8$ is the natural precursor for the weak $SU(2)$ gauge group, and 
the most plausible finite precursor to the strong $SU(3)$ gauge group is the tetrahedral group $A_4$. All these groups 
combine in the
binary tetrahedral group $2A_4$, in whose group algebra we see also the Lorentz group as an emergent symmetry group.
Of particular note is the fact that the group contains four distinct copies of $Z_3$, giving rise to four distinct generation symmetry groups,
one for each of the four types of fundamental fermions. 

I have shown how to derive
two 
more of the mixing angles that are (apparently) visible in the group algebra of the binary tetrahedral group, 
namely the mixing angle between second and third generation quarks, and the CP-violating phase 
of the CKM matrix. 
The compact part of the group algebra has exactly enough room for the five remaining mixing angles of the standard model, plus
one more that may go beyond the standard model in some way. I have suggested that the new mixing angle might
be identified as a gravi-weak mixing angle, and have suggested a way of using retarded gravity (i.e. gravity that travels at the
speed of light) to estimate its value in Earthbound experiments.
Further work will be required to find the remaining angles explicitly.

In this paper, I have examined what I believe to be essentially the unique possible mathematical model of a discrete algebraic universe,
in the hope that it will have something useful to say about the seemingly intractable problems in the foundations of physics.
I have shown how all the essential ingredients of the standard theories of both classical and quantum physics arise from this
finite model, and discussed at a general level the relationships between them. 
I have traced the conflict between quantum mechanics and
relativity to a conflict in interpretations, that does not necessarily affect the mathematics of either theory, and sketched a possible way to
resolve this conflict.

I have not done all the necessary detailed calculations to show that the proposed model reproduces the standard models exactly,
so there is still room for doubt as to whether the model I propose is viable. Nevertheless,
 I have shown that the proposed model has enough complexity to incorporate all of the
 subtleties of the standard model of particle physics, including the $24$ dimensionless parameters. 
The proposed algebraic model seems to be consistent with the experimental fact that the standard model is essentially a complete and correct
theory of everything. As for physics `beyond the standard model', the model 
explains what experiment has 
demonstrated, namely that, essentially, there is none.

Similarly, it appears to be consistent with general relativity as a theory of gravity, but suggests that the simplifying assumptions
made in practical calculations break down in extreme circumstances.
The model suggests that incorporating a contribution to gravity from the rotation of
a (not spherically symmetric) gravitating body, 
and taking into account the finite speed of propagation of gravitational waves, may be sufficient to account for the
anomalous rotation of stars in the outer regions of galaxies, that was the original reason for the hypothesis of
dark matter. If this is not enough, then the model has room to distinguish gravitational mass from inertial mass,
and to incorporate a scale factor between the two that is dependent on 
properties of the motion of the observer.

So what remains of Einstein's `castle in the air'? The finite model has, if anything, given it slightly firmer foundations. 
Special relativity has always been built on solid rock, but the foundations of general relativity are more fluid. 
I see general relativity, therefore, not as a castle in the air, but as an aeroplane,
that stays up despite having no visible means of support. My proposed model provides some support,
in the form of a neutrino wind, to keep it flying.


\begin{thebibliography}{99}
\bibitem{Einstein1935} A. Einstein (1935), Letter to Paul Langevin.
\bibitem{Einstein1954} A. Einstein (1954), Letter to Michele Besso.
\bibitem{perspective}
R. A. Wilson (2020), A group-theorist's perspective on symmetry groups in physics,
arXiv:2009.14613v5. 

 \bibitem{GG} H. Georgi and S. Glashow (1974), 
 Unity of all elementary-particle forces, {\it Physical Review Letters} {\bf 32} (8), 438.
 
\bibitem{PatiSalam} J. C. Pati and A. Salam (1974), Lepton number as the fourth `color',
{\it Phys. Rev. D} {\bf 10} (1), 275--289.

\bibitem{MDW} C. A. Manogue, T. Dray and R. A. Wilson (2022), Octions: an $E_8$ description of the standard model,
{\it J. Math. Phys.} {\bf 63}, 081703.

\bibitem{Chester} D. Chester, A. Marrani and M. Rios (2020),
Beyond the standard model with six-dimensional spacetime,
arXiv:2002.02391.

\bibitem{Lu} W. Lu (2011), Yang--Mills interactions and gravity in terms of Clifford algebra,
{\it Adv. Appl. Clifford Alg.} {\bf 21}, 145--163.

\bibitem{Stoica} O. C. Stoica (2018),  The standard model algebra---leptons, quarks and gauge from the
complex Clifford algebra $Cl(6)$, {\it Adv. Appl. Clifford Alg.} {\bf 28}, 52.

\bibitem{Clifford} R. A. Wilson (2021), 
On the problem of choosing subgroups of Clifford algebras for applications in fundamental physics,
{\it Adv. Appl. Clifford Alg.} {\bf 31}, 59.


\bibitem{MD} C. A. Manogue and T. Dray (2010), 
Octonions, E6, and particle physics, {\it  J. Phys. Conf. Ser.} {\bf 254}, 012005.

\bibitem{Todorov} I. Todorov and M. Dubois-Violette (2018), 
Deducing the symmetry of the standard model from the automorphism and structure groups of the exceptional Jordan algebra,
{\it Int. J. Mod. Phys. A} {\bf 33}, 1850118

\bibitem{GuGu} M. G\"unaydin and F. G\"ursey (1973), Quark structure and the octonions,
{\it J. Math. Phys.} {\bf 14}, 11.

\bibitem{GDixon} G. Dixon (2004), Division algebras: family replication, {\it J. Math. Phys.} {\bf 45}, 3878.

\bibitem{CFurey} C. Furey (2014), 
Generations: three prints, in colour,
{\it JHEP} {\bf 10}, 046.

\bibitem{NFurey} N. Furey (2018), Three generations, two unbroken gauge symmetries, and one 8-dimensional algebra,
{\it Phys. Lett. B} {\bf 785}, 84--89.

 \bibitem{Zee}
 A. Zee (2016), {\it Group theory in a nutshell for physicists},
 Princeton University Press.
  
\bibitem{WoitQFT} P. Woit (2017), {\it Quantum theory, groups and representations}, Springer.

\bibitem{JamesLiebeck} G. D. James and M. W. Liebeck (2012), {\it Representations and characters of groups}, 2nd ed,
Cambridge UP.

\bibitem{Pontecorvo} B. Pontecorvo (1958), Inverse beta processes and non-conservation of lepton charge,
{\it Soviet Physics JETP} {\bf 7}, 172.

\bibitem{MNS}  Z. Maki, M. Nakagawa and S. Sakata (1962),
Remarks on the unified model of elementary particles,
{\it Progress of Theoretical Physics} {\bf 28} (5), 870.

\bibitem{CODATA2018} E. Tiesinga, P. J. Mohr, D. B. Newell and B. N. Taylor (2021),
CODATA recommended values of the fundamental physical constants: 2018,
{\it Rev. Mod. Phys.} {\bf 93}, 025010.

\bibitem{WZ} 
T. Aaltonen et al. (CDF Collaboration) (2022),
High precision measurement of the $W$ boson mass with the CDF II detector, 
{\it Science} {\bf 376}:6589, 170--176.

\bibitem{SternGerlach}  W. Gerlach and O. Stern (1922), Der Experimentelle Nachweis der Richtungsquantelung im Magnetfeld,
{\it Zeitschrift f\"ur Physik}, {\bf 9}:1, 349--352.

\bibitem{grouprings} G. D. James and M. W. Liebeck (2012), {\it Representations and characters of groups}, 2nd ed,
Cambridge UP.

\bibitem{Wu} C. S. Wu, E. Ambler, R. W. Hayward, D. D. Hoppes and R. P. Hudson (1957),
Experimental test of parity conservation in beta decay,
{\it Phys. Rev.} {\bf 105} (4), 1413--1415.
 \bibitem{Griffiths} D. Griffiths (2008), {\it Introduction to elementary particles},
2nd ed, Wiley.

\bibitem{GellMann} M. Gell-Mann (1961), The eightfold way: a theory of strong interaction symmetry,
Synchrotron Lab. Report CTSL-20, Cal. Tech.

\bibitem{QCD} W. Greiner, S. Schramm and E. Stein (2007), {\it Quantum Chromodynamics}, Springer.


\bibitem{E6}  C. A. Manogue and T. Dray (2010), Octonions, $E_6$, and particle physics,
{\it J. Phys: Conf. Ser.} {\bf 254}, 012005.
\bibitem{E8} 
 A. G. Lisi (2007), An exceptionally simple theory of everything,
 arXiv:0711.0770
 

\bibitem{DG} J. Distler and S. Garibaldi (2010), There is no $E_8$ theory of everything,
{\it Communications in Math. Phys.} {\bf 298} (2), 419--436.

\bibitem{Furey} C. Furey (2014), Generations: three prints, in colour, 
{\it J. High Energy Phys.} {\bf 10}, 046.

\bibitem{DMW} T. Dray, C. A. Manogue and R. A. Wilson (2014), A symplectic representation of $E_7$,
{\it Comment. Math. Univ. Carolin.} {\bf 55}, 387--399.

\bibitem{Cabibbo} N. Cabibbo (1963),
Unitary symmetry and leptonic decays,
{\it Physical Review Letters} {\bf 10} (12), 531--533.

\bibitem{KM} M. Kobayashi and T. Maskawa (1973), CP-violation in the renormalizable theory of weak interaction,
{\it Progress of Theoretical Physics} {\bf 49} (2), 652--657.


\bibitem{tetrions} R. A. Wilson (2023), Tetrions: a discrete approach to the standard model,
arXiv:2301.11727.

\bibitem{Dirac} P. A. M. Dirac (1928), The quantum theory of the electron,
{\it Proc. Roy. Soc A} {\bf 117}, 610--624.

\bibitem{remarks} R. A. Wilson (2022),
Remarks on the group-theoretical foundations of particle physics,
{\it Intern. J. of Geom. Methods in Modern Phys.} {\bf 19}, 2250164.





\bibitem{Weinberg} S. Weinberg (1967), A model of leptons,
{\it Phys. Rev. Lett.} {\bf 19}, 1264--66.


\bibitem{SR} T. Dray (2012), {\it The geometry of special relativity}, A K Peters.



\bibitem{octahedral} R. A. Wilson (2021), Options for a finite group model of quantum mechanics,
arXiv:2104.10165.

\bibitem{gl23model} R. A. Wilson (2021), Potential applications of modular representation theory to quantum mechanics,
arXiv:2106.00550.

\bibitem{GR1} A. Einstein (1916), Die Grundlage der allgemeinen
Relativit\"atstheorie, {\it Annalen der Physik} {\bf 49} (7), 769--822.

 \bibitem{GR2} A. Einstein (1955),
{\it The meaning of relativity}, 5th ed., Princeton UP.

 \bibitem{GR3}
 G. 't Hooft (2001),
 {\it Introduction to general relativity},
 Rinton.

\bibitem{flyby} J. D. Anderson, J. K. Campbell, J. E. Ekelund, J. Ellis and J. F. Jordan (2008),
Anomalous orbital-energy changes observed during spacecraft flybys of Earth,
{\it PRL} {\bf 100}, 091102.

\bibitem{Hafele} J. C. Hafele (2013), Causal version of Newtonian theory by time-retardation
of the gravitational field explains the flyby anomalies,
{\it Progress in Physics} {\bf 9} (2), 3--8.


\bibitem{Milgrom1}  M. Milgrom (1983),
A modification of the Newtonian dynamics as a possible alternative to the
hidden mass hypothesis,
{\it Astrophysical J.} {\bf 270},  365--370.

\bibitem{Milgrom2} M. Milgrom (1983), A modification of the Newtonian dynamics:
implications for galaxies, {\it Astrophysical J.} {\bf 270}, 371--383.

\bibitem{Milgrom3} M. Milgrom (1983), A modification of the Newtonian dynamics:
implications for galaxy systems, {\it Astrophysical J.} {\bf 270}, 384--389.

\bibitem{MOND} B. Famaey and S. McGaugh (2012), Modified Newtonian dynamics (MOND):
observational phenomenology and relativistic extensions, {\it Living reviews in relativity} {\bf 15}, 10.
arXiv:1112.3960.

\bibitem{MOND2} M. Milgrom (2014), MOND laws of galactic dynamics,
{\it Monthly Notices  Roy. Astro. Soc.} {\bf 437}, 2531--41.



\bibitem{Yahalom} A. Yahalom (2019), The effect of retardation on galactic rotation curves,
{\it J. Phys.: Conf. Ser.} {\bf 1239}, 012006.

\bibitem{GL4R1} M. Leclerc (2006), The Higgs sector of gravitational gauge theories, {\it Annals of Physics} {\bf 321}, 708--743.

\bibitem{GL4R2}
D. Ivanenko and G. Sardanashvily (1983), The gauge treatment of gravity, {\it Physics Reports} {\bf 94}, 1--45.


\bibitem{oscillation}
B. Pontecorvo (1968),
Neutrino experiments and the problem of conservation of leptonic charge,
{\it Soviet Phys. JETP} {\bf 26}, 984--988.

\bibitem{neutrinos}
V. Gribov and B. Pontecorvo (1969),
Neutrino astronomy and lepton charge,
{\it Phys. Rev. D} {\bf 22} (9), 2227--2235.

\bibitem{SNO} A. Bellerive et al. (SNO collaboration) (2016),
The Sudbury neutrino observatory, {\it Nuclear Phys. B} {\bf 908}, 30--51. arXiv:1602.02469.

\bibitem{freedom1} D. J. Gross and F. Wilczek (1973), Ultraviolet behaviour of non-abelian gauge theories,
{\it PRL} {\bf 30} (26), 1343--1346.


\bibitem{freedom2} H. D. Politzer (1973), Reliable perturbative results for strong interactions,
{\it PRL} {\bf 30} (26), 1346--1349.

\bibitem{EPR} A. Einstein, B. Podolsky and N. Rosen (1935), Can quantum-mechanical description of physical reality be considered complete? 
{\it PRL} {\bf 47} (10), 777--780.

\bibitem{INI13} R. A. Wilson (2020), Possible emergence of fundamental constants from a generally covariant model of quantum mechanics,
Preprint 19013, Isaac Newton Institute, Cambridge.

\bibitem{gravimass} A. Einstein (1919), Do gravitational fields play an essential role in the structure of elementary particles of matter?
{\it Sitzungsberichte der K\"oniglichen Preu\ss ischen Akademie der Wissenschaften}, 349--356.


 \bibitem{CP}
 J. H. Christenson, J. W. Cronin, V. L. Fitch and R. Turlay (1964),
 Evidence for the $2\pi$ decay of the ${K_2}^0$ meson,
 {\it Phys. Rev. Lett.} {\bf 13}, 138.

\bibitem{ColemanGlashow} S. Coleman and S. L. Glashow (1961), {\it Phys. Rev. Lett.} {\bf 6}, 423.


\end{thebibliography}
\end{document}